\documentclass[11pt]{article}
\date{}
\makeatletter
\@addtoreset{equation}{section}
\makeatother
\usepackage{amssymb}
\usepackage{color}
\usepackage[colorlinks,linkcolor=blue,citecolor=blue]{hyperref}
\usepackage{cite}
\usepackage{graphicx}
\usepackage{amsmath}
\usepackage{times}
\usepackage{amstext}
 \usepackage[applemac]{inputenc}
\usepackage{multirow}
\usepackage{array}
\usepackage{booktabs}
\usepackage{cases}

\usepackage{pgf}
\usepackage{mathtools}
\usepackage{tikz}
\usepackage{psfrag}
\usepackage{xcolor}
\usepackage{pifont}

\usepackage{enumerate}

 \numberwithin{equation}{section}
 \newtheorem{thm}{Theorem}[section]
\newtheorem{lem}[thm]{Lemma}
\newtheorem{rem}[thm]{Remark}

\allowdisplaybreaks

\date{}

\newcommand{\diiv}{{\rm div}\ }
\newcommand{\hu}{\textbf{u}}

\newcommand{\hv}{\textbf{v}}

\newcommand{\intR}{ \int_{\mathbb{R}^3}}
\newcommand{\RN}{\mathbb{R}^N}
\newcommand{\RT}{\mathbb{R}^3}

\newcommand{\hbx}{\hfill$\Box$}
\newcommand{\pf}{{\bf Proof.\ }}

\newcommand{\X}{X_{\mathcal{SO}}}
\newcommand{\DF}{\mathcal{D}_\mathcal{F}}

\newcommand{\I}{\mathcal{I}}
\newcommand{\J}{\mathcal{J}}

\begin{document}

\title{\Large On the double critical Maxwell equations}

 \author{\small  Cong Wang\textsuperscript{1}  \ \ \ \ \ \ Jiabao Su\textsuperscript{2}\footnote{Corresponding author.\  \  E-mail addresses: sujb@cnu.edu.cn.}\\
 {\small  $^1$School of  Mathematics, Southwest Jiaotong University}\\ {\small Chengdu  610031, People's Republic of China}\\
 {\small  $^2$School of Mathematical Sciences,
 Capital Normal University}\\ {\small Beijing 100048, People's Republic of China}  }

\maketitle

\begin{abstract}
 In this paper, we focus on (no)existence  and asymptotic behavior of solutions for the double critical Maxwell equation involving with the Hardy, Hardy-Sobolev, Sobolev critical exponents. The existence and noexistence of solutions completely depend on the  power exponents and coefficients of equation. On one hand, based on the concentration-compactness ideas, applying the Nehari manifold and the mountain pass theorem, we prove the existence of the ground state solutions for the critical Maxwell equation for three different scenarios. On the other hand, for the case $\lambda<0$ and $0\leq s_2<s_1<2$, which is a type  open problem raised by Li and Lin. Draw support from a changed version of Caffarelli-Kohn-Nirenberg inequality, we find that there exists a constant $\lambda^*$ which is a negative number having explicit expression, such that the problem has no nontrivial solution as the coefficient $\lambda<\lambda^*$. Moreover, there exists  a constant $\lambda^*<\lambda^{**}<0$ such that, as $\lambda^{**}<\lambda<0$, the equation has a nontrivial solution using truncation methods. Furthermore, we establish the  asymptotic behavior of solutions of equation as coefficient converges to zero  for the all cases above.  
{\bf 2020 Mathematics Subject Classification Primary:}\ 35A15; 35B33; 35Q61.
 \end{abstract}

\section{Introduction}
In the present paper we focus on the double critical  Maxwell equation
\begin{eqnarray}
\nabla\times(\nabla \times \hu)=\frac{|\hu|^{4-2s_1}\hu}{|x|^{s_1}}+\lambda\frac{|\hu|^{4-2s_2}\hu}{|x|^{s_2}} \ \ \ {\rm in}\ \mathbb{R}^3,
\label{E}
\end{eqnarray}
where $\nabla\times(\nabla \times \cdot)$ is the curl-curl operator, $\hu:\mathbb{R}^3\to \mathbb{R}^3$ is a vector function, $0 \leq s_i\leq 2 (i=1,2)$ are constants, the number   $2^*(s):=6-2s$ is named as the Hardy (resp. Hardy-Sobolev, Sobolev) critical exponent as $s=2$ (resp. $0<s<2$, $s=0$) due to a reason that the only  continuous
embedding \begin{eqnarray}
D^{1,2}(\RT)\hookrightarrow L^{2^*(s)}(\RT;|x|^{-s}),\label{embedding}
\end{eqnarray}
which is  noncompact. The general version of the Maxwell equation is formulated as
% Faraday's Law ∑®¿≠µ⁄∂®¬…£¨ Ampere's \ Law£¨∞≤≈‡∂®¬…£ªGauss Electric Law ∏ﬂÀπµÁ≥°∂®¬…£¨ Gauss'\ Magnetic\ Law ∏ﬂÀπ¥≈≥°∂®¬…
\begin{eqnarray*}
\left\{
  \begin{array}{ll}
\partial_t \mathcal{B}+\nabla \times \mathcal{E}=0\ &{\rm (Faraday's \ Law)},\\
\nabla \times \mathcal{H}=\mathcal{J}+\partial_t \mathcal{D}\ &{\rm (Ampere's \ Law)},\\
\diiv(\mathcal{D})=\rho\ &{\rm (Gauss'\ Electric\ Law)},\\
\diiv(\mathcal{B})=0\ &{\rm (Gauss'\ Magnetic\ Law)}\\
  \end{array}
\right.
\end{eqnarray*}
with the electric field $\mathcal{E}$, the electric displacement field $\mathcal{D}$, the magnetic field $\mathcal{H}$, the magnetic induction $\mathcal{B}$, the current intensity $\mathcal{J}$ and the scalar charge density $\rho$. These fields are related by the constitutive equations determined by the material. Considering the constitutive relation $\mathcal{D}=\epsilon \mathcal{E}+\mathcal{P}_{NL}(x,\mathcal{E}),\mathcal{B}=\mu \mathcal{H}-\mathcal{M}$, where   $\epsilon=\epsilon(x)\in \mathbb{R}^{3\times3}$ is the (linear) permittivity tensor of the material,  $\mathcal{P}_{NL}$ is the nonlinear part of the polarization, $\mu=\mu(x)\in \mathbb{R}^{3\times3}$ denotes the magnetic permeability tensor and $\mathcal{M}$ the magnetization of the material.  Suppose that there are no currents, charges nor magnetization, i.e. $\mathcal{J}=0,\rho=0, \mathcal{M}=0$. Then combining with  Faraday's law and Ampere's Law, we can obtain the nonlinear electromagnetic wave equation of the form
\begin{eqnarray}
\nabla\times(\mu^{-1}\nabla \times \mathcal{E})+\epsilon(x)\partial^2_t \mathcal{E}+\partial^2_t P_{NL}(x,\mathcal{E})=0.\label{1.3}
\end{eqnarray}

The equation \eqref{1.3} is particularly challenging and in the literature there are several simplifications relying on approximation of the nonlinear electromagnetic wave equation. The most prominent one is the scalar or vector nonlinear Schr\"{o}dinger equation. In order to justify this
approximation one assumes that the term $\nabla(\diiv(\mathcal{E}))$ in $\nabla\times(\nabla\times \mathcal{E})=\nabla(\diiv(\mathcal{E}))-\Delta \mathcal{E}$ is negligible, and that one can use the so-called \emph{slowly varying envelope approximation}. % ¬˝±‰∞¸¬ÁΩ¸À∆
However, this approach may produce non-physical solutions. We can establish the time-harmonic Maxwell equation by
\begin{eqnarray}
\mathcal{E}(x,t)=\textbf{v}(x)e^{\omega t}\ {\rm for}\ x\in \mathbb{R}^3\ {\rm and}\ t\in \mathbb{R}\label{1.4}
\end{eqnarray}
with frequency $\omega>0$. We consider a special case  in \eqref{1.3} with
\begin{eqnarray*}
\mathcal{P}_{NL}(x,\mathcal{E})=-\frac{|\hv(x)|^{4-2s_1}}{|x|^{s_1}}\mathcal{E}-\gamma\frac{|\hv(x)|^{4-2s_1}}{|x|^{s_1}}\mathcal{E},  \  \mu=I,  \ \epsilon=0.
\end{eqnarray*} Then \eqref{1.3} reduces to the curl-curl equation of the type
\begin{eqnarray}
\nabla \times(\nabla \times \hv)=\omega^2\frac{|\hv|^{4-2s_1}\hv}{|x|^{s_1}}+\gamma\omega^2\frac{|\hv|^{4-2s_2}\hv}{|x|^{s_2}}\ \ {\rm in}\ \mathbb{R}^3.\label{1.5}
\end{eqnarray}
The curl operator $\nabla \times \cdot$ is challenging from the mathematical point of view and is important in mathematical physics: such an operator also  appears in the Navier-Stokes equations. Such an operator has several essential features. The kernel of $\nabla \times \cdot$ is of infinite dimensional, which makes that the corresponding energy functionals of the equation \eqref{E} or  \eqref{1.5}
\begin{eqnarray*}
J_{\omega,\gamma}(\hv)=\frac12\intR|\nabla\times \hv|^2dx-\frac{\omega^2}{6-2s_1}\intR \frac{|\hv|^{6-2s_1}}{|x|^{s_1}}dx-\frac{\gamma\omega^2}{6-2s_2}\intR \frac{|\hv|^{6-2s_2}}{|x|^{s_2}}dx
\end{eqnarray*}
is strongly indefinite, i.e. unbounded from above and below, even on the subspaces of finite codimension,   and its critical points have infinite Morse index. Besides, the Fr\'{e}chet differential of the functionals  $J_{\omega,\gamma}(\hv)$ is not sequentially weak-to weak$^*$ continuous, which leads to the limit point of a weakly convergent sequence need not to be a critical point of $J_{\omega,\gamma}(\hv)$.

It is not difficult to see that using the scaling transformation $\hu=\omega^{\frac{2}{4-2s_1}} \hv$,
 %there exists a $k_0:=\omega^{\frac{2}{4-2s_1}}$ such that as $k=k_0$,
 the equation \eqref{1.5} can be reduced to  \begin{eqnarray}
\nabla \times(\nabla \times \hu)=\frac{|\hu|^{4-2s_1}\hu}{|x|^{s_1}}+\gamma\omega^{2\frac{s_2-s_1}{2-s_1}}\frac{|\hu|^{4-2s_2}\hu}{|x|^{s_2}}\ \ {\rm in}\ \mathbb{R}^3,\label{1.6}
\end{eqnarray}
which is precisely the equation \eqref{E} with $\lambda=\gamma\omega^{2\frac{s_2-s_1}{2-s_1}}$. We denote that
\begin{eqnarray*}
J_{\lambda}(\hu)=J_{1,\lambda}(\hu).
\end{eqnarray*}

As far as we know, the initial works researching the exact solutions of the Maxwell's equation are \cite{1990Stuart,1989McleodStuartTroy}. To the best of our knowledge, the first work dealing with the Maxwell's equation using the variational methods is due to Benci and Fortunato \cite{2004Benci-Fortunato}, they introduced a series of brilliant ideas, such as splitting the function space into a divergence-free subspace and a curl-free subspace. The second attempt to tackle the Maxwell's equation is due to Azzollini, Benci, D'Aprile and Fortunato \cite{2006Azzollini-Benci-D'Aprile-Fortunato}, they used two group actions of $\mathcal{SO}:=\mathcal{SO}(2)\times\{I\}$ to  simplify the curl-curl operator $\nabla\times(\nabla\times \cdot)$ to the vector Laplacian operator $-\Delta \cdot$.  The methods used in \cite{2006Azzollini-Benci-D'Aprile-Fortunato} are to find cylindrically symmetric solutions to \eqref{1.5}.
 Let  $\mathcal{F}$ be the space of the vector fields $\hu:\mathbb{R}^3 \to \mathbb{R}^3$ such that
\begin{eqnarray}
\hu=\frac{u}{|x'|}\left(
            \begin{array}{c}
            -x_2  \\
           x_1\\
           0
            \end{array}
          \right), \  x=(x_1,x_2,x_3)\in \mathbb{R}^3 \ \ {\rm and} \ \ |x'|^2=x_1^2+x_2^2, \label{1.7}
\end{eqnarray}
where $u:\mathbb{R}^3\to \mathbb{R}$ is a $\mathcal{SO}$-invariant scalar function. Let $\mathcal{D}_\mathcal{F}:=D^{1,2}(\mathbb{R}^3,\mathbb{R}^3)\cap \mathcal{F}$.
From a direct computations we observe that  $\hu\in \mathcal{D}_\mathcal{F}$  solves the equation
\begin{eqnarray}
\nabla \times (\nabla \times \hu)=|\hu|^4\hu\ \ \ {\rm in}\ \RT,\label{1.8}
\end{eqnarray}
which is exactly the equation \eqref{E} with  $s_1=0$ and $\lambda=0$,  if and only if the function  $\phi(x):=u(|x'|,x_3)$ solves the equation
\begin{eqnarray*}
-\Delta \phi+\frac{\phi}{|x'|^2}=|\phi|^4\phi\ \ \  {\rm in}\ \mathbb{R}^3.
\end{eqnarray*}
Surprisingly, the existence of solutions of the equation \eqref{1.8} have been proved by Esteban and Lions in \cite{1989EstebanLions}. We refer to \cite{2016Bartsch-Dohnal-Plum-Reichel,2015Bartsch-Mederski,2017Bartsch-Mederski-1,2017Bartsch-Mederski-2,
2006Benci-Fortunato,2021Bieganowski,2015Mederski,2018Mederski,2016Mederski,2021MederskiSchino,2020MederskiSchinoSzulkin,2022Schino,2017Zeng} for the works on the Maxwell's equation and the references therein. Gaczkowski,   Mederski and Schino \cite{2023GaczkowskiMederskiSchino} established a ground state solution in $\mathcal{D}_\mathcal{F}$ of the equation \eqref{1.8}. Mederski and Szulkin \cite{2021MederskiSzulkin} innovatively established the optimal constants of Sobolev-type inequality for Curl operator and found the ground state solution without any symmetry assumptions.

The   equation
\begin{eqnarray}
-\Delta u=\frac{|u|^{2^*(s_1)-2}u}{|x|^{s_1}}+\lambda \frac{|u|^{2^*(s_2)-2}u}{|x|^{s_2}} \ \ {\rm in}\ \mathbb{R}^N, N \geqslant 3, \label{1.9}
\end{eqnarray}
on scalar fields has been researched widely with the  Hardy(resp. Hardy-Sobolev, Sobolev) critical exponents $2^*(s):=\frac{2(N-s)}{N-2}$ of the embedding
\begin{eqnarray}
D^{1,2}(\RN)\hookrightarrow L^{2^*(s)}(\RN;|x|^{s}),  N \geq 3 \label{embedding1}
\end{eqnarray}
as $s=2$ (resp. $0<s<2, s=0$), we refer to \cite{1976Talenti, 1976Aubin,1983Lieb,1983BN,1993ChouChu,2000GhoussoubY,2001CW,1981GS,
1997Horiuchi,2013GGN,2009Filippucci-Pucci-Robert, 2016GhoussoubRobert,2003SecchiSmetsWillem,2023WangSu-PRSE,2022WangSu-JMP,2023WangSu}. The authors \cite{2004GhoussoubKang,2010Hsia-Lin-Wadade,2012Li-Lin, 2016GhoussoubRobert, 2006GhoussoubRobert}  considered   \eqref{1.9}  in the half space $\mathbb{R}^N_+$.
In \cite{2012Li-Lin}, Li and Lin gave an open problem about \eqref{1.9} in the half space $\mathbb{R}^N_+$, which has not been fully resolved yet and which is the spindle of problems what  has been studied  in \cite{2023WangSu} and we partially answered this question.

The authors \cite{2002BadialeTarantello,2007TertikasTintarev,2008Musina,2013SunZhang,2013CarriaoDemarqueMiyagaki,2013Sun,2011Mazya} studied the equation
\begin{eqnarray}
-\Delta u=\frac{|u|^{2^*(s_1)-2}u}{|\bar{x}|^{s_1}}+\lambda \frac{|u|^{2^*(s_2)-2}u}{|\bar{x}|^{s_2}} \ \ {\rm in}\ \mathbb{R}^N, \ N \geq 3  \label{1.11}
\end{eqnarray}
with the so-called Hardy-Maz'ya(resp. Hardy-Sobolev-Maz'ya, Sobolev) critical exponent  $2^*(s):=\frac{2(N-s )}{N-2} $ of the embedding
\begin{eqnarray*}
D^{1,2}(\RN)\hookrightarrow L^{2^*(s)}(\RN;|\bar{x}|^{s}),
\end{eqnarray*}
as $s=2$(resp. $0<s<2$, $s=0$), where $x=(\bar{x}, \tilde{x})\in \mathbb{R}^k\times \mathbb{R}^{N-k}, 1\leq k\leq N$.

For the equation \eqref{E}, based on the transformation \eqref{1.7}, a direct computations we observe that  $\hu\in \mathcal{D}_\mathcal{F}$  solves \eqref{E} if and only if $u\in \X$ solves
\begin{eqnarray}
-\Delta u+\frac{u}{|x'|^2}=\frac{|u|^{4-2s_1} u}{|x|^{s_1}}+\lambda \frac{|u|^{4-2s_2} u}{|x|^{s_2}}\ \ \  {\rm in}\ \mathbb{R}^3, \label{1.12}
\end{eqnarray}
\begin{eqnarray}
x:=(x',x_3)\in \mathbb{R}^2\times \mathbb{R}, \ x':=(x_1,x_2)\in \mathbb{R}^2, \ x_3\in \mathbb{R}, \label{1.13}
\end{eqnarray}
where $X_{\mathcal{SO}}$ is the closed subspace of the function space  $X$ defined by
\begin{equation*}
X=\left\{u\in D^{1,2}(\RT) \ \Big| \ \int_{\RN}\frac{|u|^2}{|x'|^2}dx<\infty\right\}
\end{equation*}
consisting of the functions invariant under the usual group action of $\mathcal{SO}:=\mathcal{O}(2)\times \{I\}\subset \mathcal{O}(3)$.
 Note that this is equivalent requiring that such functions be invariant under the action of $\mathcal{O}(2)\times \{I\}$ because for every $\xi_1,\xi_2\in \mathbb{S}^{1}$, there exists $g\in \mathcal{SO}(2)$ such that $\xi_2=g\xi_1$, where $\mathcal{SO}(2)\subset \mathcal{O}(2)$ stands for the special orthogonal group in $\mathbb{R}^2$.
 The space $X$ is a Hilbert space endowed with the scalar product \begin{eqnarray*}
\langle u,v\rangle\in X\times X\mapsto \int_{\RT}\nabla u\cdot\nabla v+\frac{uv}{|x'|^2}dx
\end{eqnarray*}
and the corresponding norm
\begin{eqnarray*}
\|u\|:=\langle u,u\rangle^{1/2}.
\end{eqnarray*}
We remark that, for any $\hu\in  \mathcal{D}_\mathcal{F} $, on the sense of \eqref{1.7} with $u\in \X$,  there holds that \begin{eqnarray}
\left.
  \begin{array}{ll}
\displaystyle \int_{\RT}|\nabla\times\hu|^2dx=\int_{\RT}|\nabla u|^2+\frac{|u|^2}{|x'|^2}dx=:\mathcal{A}(u),\\ [1em]
\displaystyle \int_{\RT}\frac{|\hu|^{6-2s_1}}{|x|^{s_1}}dx=\int_{\RT}\frac{|u|^{6-2s_1}}{|x|^{s_1}}dx=:\mathcal{B}(u),\\ [1em]
\displaystyle  \int_{\RT}\frac{|\hu|^{6-2s_2}}{|x|^{s_2}}dx=\int_{\RT}\frac{|u|^{6-2s_2}}{|x|^{s_2}}dx=:\mathcal{C}(u).
  \end{array}
\right.\label{1.14}
\end{eqnarray}
We define the corresponding functional of the equation \eqref{1.12} as
\begin{eqnarray*}
\I_\lambda(u)=\frac{1}{2}\intR|\nabla u|^2+\frac{|u|^2}{|x'|^2}dx-\frac{1}{6-2s_1}\intR\frac{|u|^{6-2s_1}}{|x|^{s_1}}dx-\frac{\lambda}{6-2s_2}\intR\frac{|u|^{6-2s_2}}{|x|^{s_2}}dx
\end{eqnarray*}
and define the    corresponding Nehari manifold
\begin{eqnarray*}
\mathcal{N}:=\left\{u\in \X \ | \ \langle \I'_\lambda(u),u \rangle=0\right\}.
\end{eqnarray*}
According \eqref{1.14}, we see that  $J_{\lambda}(\hu)=\I_\lambda(u)$ under the transformation \eqref{1.7}, for more information and related derivation process one sees  \cite{2023GaczkowskiMederskiSchino}. {\it Thus, we will achieve the solvability of the equation \eqref{E} in $\mathcal{D}_\mathcal{F}$  by researching the equation \eqref{1.12} in $\X$}. The ground state solutions of the equation \eqref{1.12} are the extremal functions of
\begin{eqnarray*}
m_\lambda:=\inf_{u\in \mathcal{N}}\I_\lambda(u).
\end{eqnarray*}
 We define a number as
\begin{eqnarray}
\bar{\lambda}:=\inf_{u\in \X}\frac{\intR|\nabla u|^2+\frac{|u|^2}{|x'|^2}dx}{\intR\frac{|u|^2}{|x|^2}dx},\label{1.15}
\end{eqnarray}
\begin{rem}\label{rem0}
It follows from the Hardy inequality that  $\bar{\lambda}\geq \frac14$. And it is open that whether the best constant $\bar{\lambda}$ is achieve.
\end{rem}

For the  critical equation  with one Sobolev critical exponent, the case that $\lambda=0$  and  $s=0$, the existence of ground state solutions of \eqref{E} in $\DF$ has been obtained in \cite{2023GaczkowskiMederskiSchino}. Applying the quotient methods, we obtain a result as follows for  the critical equation  with one Hardy  critical exponent and one  Hardy-Sobolev critical  exponent, the  case that $0<s_1<2, $ $s_2=2$,
\begin{thm}\label{thm1.1}
Assume that $\lambda<\bar{\lambda},$ $ s_2=2$ and $0<s_1<2$. Then the equation \eqref{E} has a nontrivial ground state solution in $\DF$.
\end{thm}
For each $\lambda<\bar{\lambda}$ fixed, \begin{eqnarray*}
\|u\|_\lambda:=\left(\intR|\nabla u|^2+\frac{|u|^2}{|x'|^2}-\lambda \frac{|u|^2}{|x|^2}dx\right)^{1/2}
\end{eqnarray*} defines an equivalent norm to $\|u\|$ on $\X$.

For the double critical equations, the case $\lambda\not=0$, we  obtain  the  following results in the cases that $\lambda>0$ and $\lambda<0$.
\begin{thm}\label{thm1.2} Assume that  $\lambda>0$, $0\leq s_1<s_2<2$. Then the equation \eqref{E} has a nontrivial ground state solution in $\DF$.
\end{thm}

\begin{thm}\label{thm1.3}
Assume that $\lambda<0$, $0<s_1<s_2<2$. Then the equation \eqref{E} has a nontrivial ground state solution in $\DF$.
\end{thm}

Define  a real number
\begin{eqnarray*}
\lambda^*:=-\frac{2-s_1}{2-s_2}\bar{S}^{6-2s_1}\left(\frac{s_1-s_2}{2-s_2} \ \bar{S}^{6-2s_1}\right)^{\frac{(2-s_1)(s_1-s_2)}{(2-s_2)^2}},
\end{eqnarray*}
where $\bar{S}$ is a best constant of Caffarelli-Kohn-Nirenberg type inequality, see Lemma \ref{lemCKN} in Section 3.
 We have the following nonexistence  and existence results of the nontrivial solutions of \eqref{E} in the case $\lambda<0$.

\begin{thm}\label{thm1.4} Assume that $\lambda<\lambda^*$, $0\leq s_2<s_1<2$. Then the equation \eqref{E} has only zero solution in $\DF$.
\end{thm}
\begin{thm}\label{thm1.5}
Assume that $0\leq s_2<s_1<2$. Then there exists a $\lambda^{**}\in(\lambda^*,0)$ such that the equation \eqref{E} has a nontrivial solution in $\DF$ as  $\lambda^{**}<\lambda<0$.
\end{thm}

It follows from results above, we summarize the  solvability of the equation \eqref{E} as follows.
\begin{center}
{\footnotesize
\begin{tabular}{ccc}
\toprule[1mm]
$\lambda$ &  $s_1, s_2$ & Does the equation \eqref{E} have nontrivial solutions in $\DF$ \\
\midrule[1pt] \specialrule{0em}{2pt}{2pt}
$\lambda=0$ & $s_1=0$ &Yes(see \cite{2023GaczkowskiMederskiSchino})\\
\hline \specialrule{0em}{2pt}{2pt}
$\lambda<\bar{\lambda}$ & $0<s_1<2$, $s_2=2$ & Yes(see Theorem \ref{thm1.1}) \\
\hline \specialrule{0em}{2pt}{2pt}
$\lambda>0$& $0\leq s_1<s_2<2$ & Yes(see Theorem \ref{thm1.2})\\
\hline \specialrule{0em}{2pt}{2pt}
$\lambda<0$&  $0<s_1<s_2<2$ & Yes(see Theorem \ref{thm1.3}) \\
 \hline \specialrule{0em}{2pt}{2pt}
$\lambda<\lambda^*$ & $0\leq s_2<s_1<2$ & No(see Theorem \ref{thm1.4}) \\
\hline \specialrule{0em}{2pt}{2pt}
$\lambda^{**}<\lambda<0$& $0\leq s_2<s_1<2$ & Yes(see Theorem \ref{thm1.5}) \\
\bottomrule[1mm]
\end{tabular}}
\end{center}
\begin{rem}\label{rem1}
According to the table above, it is not difficulty to find that the solvability of the equation \eqref{E} in $\DF$ remains unsolved under the following cases,
\begin{itemize}
 \item  $\lambda\geq \bar{\lambda}, 0\leq s_1<2, s_2=2$,
  \item $\lambda<\bar{\lambda}$ and $\lambda\not=0$, $s_1=0, s_2=2$,
  \item $\lambda<0, 0=s_1<s_2<2$.
\end{itemize}
\end{rem}

Now, we investigate the  asymptotic behavior  of the solutions of \eqref{E} as $\lambda\to 0$.
\begin{thm}\label{thm1.6}  Assume that  $0<s_1<2$ and $s_2=2$.  Then there exist a sequence $\{\lambda_n<\bar{\lambda}\}$ and $\hu\in \DF$, a ground state solution  of   \eqref{E}  in $\DF$ with $\lambda=0$,  such that    the solution sequence  $\{\hu_n\} \subset \DF$ of    \eqref{E} corresponding to the sequence $\{\lambda_n<\bar{\lambda}\}$, satisfies  that $\hu_n\to \hu$  in $\DF$ as  $\lambda_n\to0$. \end{thm}

\begin{thm}\label{thm1.7}
Assume  that $0<s_1<s_2<2, \lambda>0$. Then there exist a sequence $\{\lambda_n>0\}$ and $\hu\in \DF$, a ground state solution  of   \eqref{E}  in $\DF$ with $\lambda=0$,  such that    the solution sequence  $\{\hu_n\} \subset \DF$ of    \eqref{E} corresponding to the sequence $\{\lambda_n>0\}$, satisfies  that $\hu_n\to \hu$  in $\DF$ as  $\lambda_n\to0^+$.
\end{thm}

\begin{thm}\label{thm1.8}
Assume  that $0<s_1<s_2<2, \lambda<0$.  Then there exist a sequence $\{\lambda_n<0\}$ and $\hu\in \DF$, a ground state solution  of   \eqref{E}  in $\DF$ with $\lambda=0$,  such that    the solution sequence  $\{\hu_n\} \subset \DF$ of    \eqref{E} corresponding to the sequence $\{\lambda_n<0\}$, satisfies  that $\hu_n\to \hu$  in $\DF$ as  $\lambda_n\to0^-$.
\end{thm}

\begin{thm}\label{thm1.9}
  Assume that $0\leq s_2<s_1<2, \lambda<0$. Then there exist a sequence $\{\lambda_n<0\}$ and $\hu\in \DF$, a ground state solution  of   \eqref{E}  in $\DF$ with $\lambda=0$,  such that    the solution sequence  $\{\hu_n\} \subset \DF$ of    \eqref{E} corresponding to the sequence $\{\lambda_n<0\}$, satisfies  that $\hu_n\to \hu$  in $\DF$ as  $\lambda_n\to0^-$.
 \end{thm}

\begin{rem}\label{rem}
As arguments above, let $\mathcal{F}$ be the space of the vector fields $\hu:\mathbb{R}^3\to \mathbb{R}^3$ such that
\begin{eqnarray}
\hu=\frac{u}{|x'|}\left(
            \begin{array}{c}
            -x_2  \\
           x_1\\
           0
            \end{array}
          \right), x=(x_1,x_2,x_3)\in \mathbb{R}^3\ {\rm and} \ |x'|^2=x_1^2+x_2^2, \label{1.16}
\end{eqnarray}
where $u\in\X$  and  $\mathcal{D}_\mathcal{F}:=D^{1,2}(\mathbb{R}^3,\mathbb{R}^3)\cap \mathcal{F}$, from a direct computations we observe that  $\hu\in \mathcal{D}_\mathcal{F}$  solves \eqref{E} if and only if $u\in \X$ solves
\begin{eqnarray}
-\Delta u+\frac{u}{|x'|^2}=\frac{|u|^{4-2s_1}u}{|x|^{s_1}}+\lambda \frac{|u|^{4-2s_2}u}{|x|^{s_2}} \ \ {\rm in}\ \mathbb{R}^3. \label{1.17}
\end{eqnarray}
Thus, the proofs  of Theorems {\rm\ref{thm1.1}--\ref{thm1.9}} depend on the corresponding results about the equation \eqref{1.17}.
 For more  information one  can refer to the important paper {\rm\cite{2023GaczkowskiMederskiSchino}}.
\end{rem}

\subsection*{Furthermore problems}
In this subsection, we give two related problems, based on the proceed of this paper, we can get the similar results of Theorems \ref{thm1.1}--\ref{thm1.9}.

Let $N\geq 3, k\geq2$, we denote
\begin{eqnarray}
x:=(\bar{x},\tilde{x})\in \mathbb{R}^k\times \mathbb{R}^{N-k}, \bar{x}=(x_1,x_2,\cdots, x_k)\in \mathbb{R}^k, \tilde{x}=(x_{k+1}, \cdots, x_N)\in \mathbb{R}^{N-k}.\label{1.18}
\end{eqnarray}

The first problem is the following

\subsubsection*{$(\mathrm{\textbf{P}}_1)$ The double critical Maxwell equation in higher dimensions.}
 Consider the case $k=2$ in \eqref{1.18}.  In order to find a suitable counterpart for the curl-curl operator $\nabla\times(\nabla\times \cdot)$ in higher dimensions, we can use the identity
\begin{eqnarray*}
\nabla\times(\nabla\times \hu)=\nabla(\nabla\cdot\hu)-\Delta \hu,  \ \ \hu\in C^2(\RN; \RN)
\end{eqnarray*}
to research the equation
\begin{eqnarray}
\nabla\times(\nabla\times \hu)=\frac{|\hu|^{2^*(s_1)-2}\hu}{|x|^{s_1}}+\lambda\frac{|\hu|^{2^*(s_2)-2}\hu}{|x|^{s_2}}\ \ {\rm in}\ \mathbb{R}^N, \label{1.19}
\end{eqnarray}
where $2^*(s):=\frac{2(N-s)}{N-2}$ is the critical exponent of the embedding \eqref{embedding1}. We can find a solution of \eqref{1.19} in $\mathcal{D}_\mathcal{F}$, where $\mathcal{D}_\mathcal{F}:=D^{1,2}(\mathbb{R}^N,\mathbb{R}^N)\cap \mathcal{F}$, $\mathcal{F}$ is the space of the vector fields $\hu:\mathbb{R}^N\to \mathbb{R}^N$ such that
\begin{eqnarray}
\hu=\frac{u}{r}\left(
            \begin{array}{c}
            -x_2  \\
           x_1\\
           0
            \end{array}
          \right), x=(x_1,x_2,\tilde{x})\in \mathbb{R}^N\ {\rm and} \ r^2=x_1^2+x_2^2\label{1.20}
\end{eqnarray}
for some $\mathcal{SO}:=\mathcal{O}(2)\times \{I_{N-2}\}-$invariant scalar function $u:\mathbb{R}^N\to \mathbb{R}$ and $u\in X_{\mathcal{SO}}$, where $\X$ is the subspace of
\begin{eqnarray*}
X:=\left\{u\in D^{1,2}(\RN)\Big|\ \int_{\RN}\frac{|u|^2}{r^2}dx<\infty\right\}
\end{eqnarray*}
consisting of the functions invariant under the usual action of $\mathcal{SO}$. Then,  $\hu\in \mathcal{D}_\mathcal{F}$  solves \eqref{1.19} if and only if $u\in\X$ solves
\begin{eqnarray}
-\Delta u+\frac{u}{r^2}=\frac{|u|^{2^*(s_1)-2} u}{|x|^{s_1}}+\lambda \frac{|u|^{2^*(s_2)-2} u}{|x|^{s_2}}\ \ \  {\rm in}\ \mathbb{R}^N.\label{1.21}
\end{eqnarray}
It is worth mentioning that Schino \cite[Corollary 4.1.4]{2022Schino}  researched \eqref{1.19} with $s_1=0$ and $\lambda=0$.

The second problem is the following
\subsubsection*{$(\mathrm{\textbf{P}}_2)$ The double critical semilinear equation in higher dimensions}
Consider the double critical equation
\begin{eqnarray}
-\Delta u+\frac{u}{|\bar{x}|^2}=\frac{|u|^{2^*(s_1)-2} u}{|x|^{s_1}}+\lambda \frac{|u|^{2^*(s_2)-2} u}{|x|^{s_2}}\ \ \  {\rm in}\ \mathbb{R}^N.\label{1.22}
\end{eqnarray}
As the case $k=2$ in \eqref{1.18}, the equation \eqref{1.22} turns into the equation \eqref{1.21}. Indeed, based on the proceed of the present paper,  we are able to verify the similar results of Theorem \ref{thm1.1}-\ref{thm1.9} for the equation \eqref{1.22} in the following function space $X$, without the  restriction that the functions are invariant under the usual action of $\mathcal{SO}:=\mathcal{O}(2)\times \{I_{N-2}\}$,
\begin{equation*}
X=\left\{u\in D^{1,2}(\RN) \ \Big|\ \int_{\RN}\frac{|u|^2}{|\bar{x}|^2}dx<\infty\right\}.
\end{equation*}
Since the quantity $\frac{|\phi|^2}{|\bar{x}|^2}$ need not be integrable for $\phi\in C_0^\infty(\RN)$, $C_0^\infty(\RN)\not\subset X$.
 Thus  we can not find a solution of the equation \eqref{1.22} in $D^{1,2}(\RN)$ in the methods in the present paper.

 For the case $k>2$ in \eqref{1.18}.  On one hand, we also prove the similar results of Theorems \ref{thm1.1}--\ref{thm1.9} for the equation \eqref{1.22} in the function space
\begin{equation*}
X=\left\{u\in D^{1,2}(\RN)\Big|\ \int_{\RN}\frac{|u|^2}{|\bar{x}|^2}dx<\infty\right\}.
\end{equation*}
On the other hand, since, see \cite{2002BadialeTarantello},
\begin{eqnarray*}
\int_{\RN}\frac{|u|^2}{|\bar{x}|^2}dx\leq \left(\frac{2}{k-2}\right)^2\int_{\RN}|\nabla u|^2dx,\ \ \ \forall \ u\in D^{1,2}(\RN),
\end{eqnarray*}
the norm
\begin{eqnarray*}
\left(\int_{\RT}|\nabla u|^2+\frac{|u|^2}{|\bar{x}|^2}dx\right)^{\frac12}
\end{eqnarray*}
of the function space
\begin{equation*}
X=\left\{u\in D^{1,2}(\RN)\Big|\ \int_{\RN}\frac{|u|^2}{|\bar{x}|^2}dx<\infty\right\}.
\end{equation*}
and the norm $\left(\int_{\RT}|\nabla u|^2dx\right)^{\frac12}$ of $D^{1,2}(\RN)$ are equivalent. Therefore, for $k\geq3$,  we can prove the existence of solutions for the equation \eqref{1.22} in $D^{1,2}(\RN)$.

Finally, we remark that, for the case  $\lambda=0$, $N\geq3, k\geq2$ with $ s:=s_1$, the equation \eqref{1.22}  becomes
\begin{eqnarray}
-\Delta u+\frac{u}{|\bar{x}|^2}=\frac{|u|^{2^*(s)-2} u}{|x|^s}\ \ \  {\rm in}\ \mathbb{R}^N.\label{1.23}
\end{eqnarray}
The similar result of Theorem \ref{thm1.1} for the equation \eqref{1.23} is new, to the best of our knowledge, none consider the case that different effective dimensions of weight functions, that is, the effective dimensions of the terms $\frac{1}{|\bar{x}|^2}$ and $\frac{1}{|x|^s}$ may be different.

\subsection*{The structure of the paper}

In Section \ref{sec2}, we prove Theorem \ref{thm1.1} by applying the quotient methods and concentration compactness ideas. In Section \ref{sec3}, the existence of the ground state solutions of \eqref{E}, that is the proofs of Theorems \ref{thm1.2} and \ref{thm1.3} is confirmed. In Section \ref{sec4}, we focus on the similar open problem raised by Li and Lin \cite{2012Li-Lin}, we prove the nonexistence (see subsection \ref{sec4.1}) and existence (see subsection \ref{sec4.2}) of nontrivial solutions contained in Theorems \ref{thm1.4} and \ref{thm1.5}.  In the final section, we establish the asymptotic behavior of solutions of  \eqref{E}, which is the proofs of Theorems \ref{thm1.6}, \ref{thm1.7}, \ref{thm1.8} and \ref{thm1.9}.

Now we give some notations  description.
\begin{itemize}
  \item Set $B_R(0)$ is a ball with center $0\in \RN$ and radius $R$ in $\RN$, specifically, as $N=1$ and $0\in \mathbb{R}$, $B_R(0):=(-R,R)$.
  \item According to the markings \eqref{1.13},  $0:=(0_1,0_2,0_3)\in \RT$ and $0':=(0_1,0_2)\in \mathbb{R}^2$.
  \item Set $b>a>0$,
\begin{itemize}
    \item    $B_{a,b}(0):=B_b(0')\times B_b(0_3) \setminus  B_a(0')\times B_a(0_3)$.
  \item $B_{a,b}(0'):=B_b(0')\setminus B_a(0')$.
  \item  $B_{a,b}(0_3):=B_b(0_3)\setminus B_a(0_3)$.
\end{itemize}
\end{itemize}

\section{Proof of Theorem \ref{thm1.1}}\label{sec2}
In this section, we focus on the proof of Theorem \ref{thm1.1}, that is the equation \eqref{E} with $\lambda<\bar{\lambda}$, where $\bar{\lambda}$ is in \eqref{1.15}, to simplify notation we write $s$ in place of $s_1$, we consider the equation
\begin{eqnarray}
\nabla\times(\nabla \times \hu)-\lambda \frac{\hu}{|x|^2}=\frac{|\hu|^{4-2s}\hu}{|x|^{s}}\ \ {\rm in}\ \mathbb{R}^3.
\label{2.1}
\end{eqnarray}
The existence of ground state solution of \eqref{2.1} in $\DF$ has been obtained in \cite{2023GaczkowskiMederskiSchino} as $s=0, \lambda=0$. Next, we will prove the result for the case $0<s<2, \lambda<\bar{\lambda}$.

{\bf Proof of Theorem \ref{thm1.1}}\  We observe that   $\hu\in\DF$  solves \eqref{2.1} if and only if $u\in \X$ solves
\begin{eqnarray*}
-\Delta \phi+\frac{u}{|x'|^2}-\lambda \frac{u}{|x|^2}=\frac{|u|^{4-2s}u}{|x|^s}\ \ \  {\rm in}\ \mathbb{R}^3.
\end{eqnarray*}
We define
\begin{eqnarray}
S_{\lambda,s}(\RT)=\inf_{u\in \X\setminus\{0\}}\frac{\int_{\RT}|\nabla u|^2+\frac{|u|^2}{|x'|^2}-\lambda \frac{|u|^2}{|x|^2}dx}{\left(\int_{\RT}\frac{|u|^{6-2s}}{|x|^s}\right)^{\frac{1}{3-s}}}.\label{Inequality}
\end{eqnarray}
 Let $\{\tilde{u}_n\}\subset \X$ be a minimizing sequence for $S_{\lambda,s}(\RT)$ such that
\begin{eqnarray*}
\int_{\RT}\frac{|\tilde{u}_n|^{6-2s}}{|x|^s}dx=1,\ \ \ \  \lim_{n\to\infty}\int_{\RT}|\nabla \tilde{u}_n|^2+\frac{|\tilde{u}_n|^2}{|x'|^2}-\lambda \frac{|u_n|^2}{|x|^2}dx=S_{\lambda,s}(\RT).
\end{eqnarray*}
For any $n$, there exists $r_n>0$ such that
\begin{eqnarray*}
\int_{B_{r_n}(0')\times B_{r_n}(0_3)} \frac{|\tilde{u}_n|^{6-2s}}{|x|^s}dx=\frac12.
\end{eqnarray*}
Define $u_n(x)=r_n^{\frac{1}{2}}\tilde{u}_n(r_nx)$, then $u_n\in \X$, and we have
\begin{eqnarray}
\lim_{n\to\infty}\int_{\RT}|\nabla u_n|^2+\frac{|u_n|^2}{|x'|^2}-\lambda \frac{|u_n|^2}{|x|^2}dx=S_{\lambda,s}(\RT).\label{2.3}
\end{eqnarray}
\begin{eqnarray}
\lim_{n\to\infty}\int_{\RT}|\nabla u_n|^2+\frac{|u_n|^2}{|x'|^2}-\lambda \frac{|u_n|^2}{|x|^2}dx=1,\ \ \ \int_{B_1(0')\times B_1(0_3)}\frac{|u_n|^{6-2s}}{|x|^s}dx=\frac{1}{2}.\label{2.4}
\end{eqnarray}
We first claim that, up to a subsequence,
\begin{eqnarray}
\lim_{R\to\infty}\lim_{n\to\infty}\int_{B_R(0')\times B_R(0_3)}\frac{|u_n|^{6-2s}}{|x|^s}dx=1, \label{2.5}
\end{eqnarray}
Indeed, for $n\in N$ and $r>0$, we define
\begin{eqnarray*}
Q_n(r):=\int_{B_r(0')\times B_r(0_3)}\frac{|u_n|^{6-2s}}{|x|^s}dx.
\end{eqnarray*}
Since $0\leq Q_n\leq1$ and $r\mapsto Q_n(r)$ is nondecreasing for all $n\in N$, then up to a subsequence, there exists $Q:[0,+\infty)\to \mathbb{R}$ nondecreasing such that $Q_n(r)\to Q(r)$ as $n\to+\infty$ for a.e. $r>0$. Set
\begin{eqnarray*}
\alpha=\lim_{r\to\infty}Q(r).
\end{eqnarray*}
It follows from \eqref{2.3} and \eqref{2.4} that $1/2\leq \alpha\leq 1$. Up to taking another subsequence, there exist $\{r_n\}_n, \{\bar{r}_n\}_n \subset (0,+\infty)$ satisfying
\begin{eqnarray*}
\left\{
  \begin{array}{ll}
\displaystyle 2r_n\leq \bar{r}_n\leq 3r_n\ \ {\rm for\ any}\ n\in \mathbb{N}^+,\\
\displaystyle  \lim_{n\to\infty}r_n=\lim_{k\to\infty}\bar{r}_n=+\infty,\\
\displaystyle  \lim_{n\to\infty}Q_n(r_n)=\lim_{n\to\infty}Q_n(\bar{r}_n)=\alpha.
  \end{array}
\right.
\end{eqnarray*}
In particular,
\begin{eqnarray}
\lim_{n\to\infty}\int_{B_{r_n}(0')\times B_{r_n}(0_3)}\frac{|u_n|^{6-2s}}{|x|^s}dx=\alpha\ {\rm and}\ \lim_{n\to\infty}\int_{\RT\setminus B_{\bar{r}_n}(0')\times B_{\bar{r}_n}(0_3)}\frac{|u_n|^{6-2s}}{|x|^s}dx=1-\alpha.\label{2.6}
\end{eqnarray}
We claim that
\begin{eqnarray}
\lim_{n\to\infty}r_n^{-2}\int_{B_{r_n}(0')\times B_{r_n}(0_3)} u_n^2dx=0.\label{2.7}
\end{eqnarray}
Indeed, for all $x\in B_{r_n,\bar{r}_n}(0)$, we have $r_n\leq |x'|\leq 2r_n$. Therefore, H\"{o}lder's inequality yieds
\begin{eqnarray*}
\int_{B_{r_n,\bar{r}_n}(0)}u_n^2dx\leq C r_n^2\left(\int_{B_{r_n,\bar{r}_n}(0)}\frac{|u_n|^{6-2s}}{|x|^s}dx\right)^{\frac{1}{3-s}}
\end{eqnarray*}
for all $n\in N$, conclusion \eqref{2.7} then follows from \eqref{2.6}.

We now let $\varphi_1\in C_0^\infty(\mathbb{R}^2), \varphi_2\in C_0^\infty(\mathbb{R})$ and
\begin{eqnarray*}
\varphi_1(x):=\left\{
                \begin{array}{ll}
         1\ \ {\rm for}\ \ x\in B_1(0'),\\
         0\ \ {\rm for}\ \ x\in \mathbb{R}^2\setminus B_2(0').
                \end{array}
              \right.
\varphi_2(x_3):=\left\{
                \begin{array}{ll}
         1\ \ {\rm for}\ \ x_3\in B_1(0_3),\\
         0\ \ {\rm for}\ \ x_3\in \mathbb{R}\setminus B_2(0_3).
                \end{array}
              \right.
\end{eqnarray*}
For $n\in N$, we define $\varphi_n(x)=\varphi_{1n}(x')\varphi_{2n}(x_3)$, where
\begin{eqnarray*}
&&\varphi_{1n}(x'):=\varphi_1\left(\frac{|x'|}{\bar{r}_n-r_n}+\frac{\bar{r}_n-2r_n}{\bar{r}_n-r_n}\right)\ \ {\rm for}\ \  x'\in \mathbb{R}^2,\\
&&\varphi_{2n}(x_3):=\varphi_2\left(\frac{|x_3|}{\bar{r}_n-r_n}+\frac{\bar{r}_n-2r_n}{\bar{r}_n-r_n}\right)\ \ {\rm for}\ \  x_3\in \mathbb{R}.
\end{eqnarray*}
One can easily check that $\varphi_{n}u_n, (1-\varphi_{n})u_n\in \X$. It follows that
\begin{eqnarray*}
\int_{\RT}\frac{|\varphi_{n}u_n|^{6-2s}}{|x|^s}dx\geq \int_{B_{r_n}(0')\times B_{r_n}(0_3)}\frac{|u_n|^{6-2s}}{|x|^s}dx=\alpha+o(1),
\end{eqnarray*}
\begin{eqnarray*}
\int_{\RT}\frac{|(1-\varphi_{n})u_n|^{6-2s}}{|x|^s}dx\geq \int_{\RT\setminus B_{\bar{r}_n}(0')\times B_{\bar{r}_n}(0_3)}\frac{|u_n|^{6-2s}}{|x|^s}dx=1-\alpha+o(1)
\end{eqnarray*}
as $n\to\infty$. The Hardy-Sobolev inequality and \eqref{2.7} implies that, as $n\to\infty$,
\begin{eqnarray*}
&&S_{\lambda,s}(\RT)\left(\int_{\RN}\frac{|\varphi_{n}u_n|^{6-2s}}{|x|^s}dx\right)^{\frac{1}{3-s}}\\
&\leq& \int_{\RT}|\nabla(\varphi_{n}u_n)|^2+\frac{|\varphi_{n}u_n|^2}{|x'|^2}-\lambda \frac{|\varphi_{n}u_n|^2}{|x|^2}dx\\
&\leq& \int_{\RT}\varphi_{n}^2\left(|\nabla u_n|^2+\frac{|u_n|^2}{|x'|^2}-\lambda \frac{|u_n|^2}{|x|^2}\right)dx+o(1).
\end{eqnarray*}
Similarly,
\begin{eqnarray*}
&&S_{\lambda,s}(\RT)\left(\int_{\RT} \frac{|(1-\varphi_{n})u_n|^{6-2s}}{|x|^s}dx\right)^{\frac{1}{3-s}}\\
&\leq& \int_{\RT}(1-\varphi_{n})^2\left(|\nabla u_n|^2+\frac{|u_n|^2}{|x'|^2}-\lambda \frac{|u_n|^2}{|x|^2}\right)dx+o(1)
\end{eqnarray*}
as $n\to\infty$. Therefore, we have that
\begin{eqnarray*}
\int_{\RT}\frac{|\varphi_{n}u_n|^{6-2s}}{|x|^s}dx\geq \int_{B_{r_n}(0')\times B_{r_n}(0_3)} \frac{|u_n|^{6-2s}}{|x|^s}dx=\alpha+o(1).
\end{eqnarray*}
\begin{eqnarray*}
\int_{\RT}\frac{|(1-\varphi_{n})u_n|^{6-2s}}{|x|^s}dx\geq \int_{\RT\setminus B_{\bar{r}_n}(0')\times B_{\bar{r}_n}(0_3)} \frac{|u_n|^{6-2s}}{|x|^s}dx=1-\alpha+o(1).
\end{eqnarray*}
To sum up we know,
\begin{eqnarray*}
&&S_{\lambda,s}(\RT)\left(\alpha^{\frac{1}{3-s}}+(1-\alpha)^{\frac{1}{3-s}}+o(1)\right)\\
&\leq& S_{\lambda,s}(\RT)\left(\left(\int_{\RT}\frac{|\varphi_{n}u_n|^{6-2s}}{|x|^s}dx\right)^{\frac{1}{3-s}}+\left(\int_{\RT} \frac{|(1-\varphi_{n})u_n|^{6-2s}}{|x|^s}dx\right)^{\frac{1}{3-s}}\right)\\
&\leq& \int_{\RT}\left(\varphi_{n}^2+(1-\varphi_{n})^2\right)\left(|\nabla u_n|^2+\frac{|u_n|^2}{|x'|^2}-\lambda \frac{|u_n|^2}{|x|^2}\right)dx+o(1)\\
&=& \int_{\RT}\left(1-2\varphi_{n}(1-\varphi_{n})\right)\left(|\nabla u_n|^2+\frac{|u_n|^2}{|x'|^2}-\lambda \frac{|u_n|^2}{|x|^2}\right)dx+o(1)\\
&\leq& S_{\lambda,s}(\RT)+2|\lambda |\int_{\RT}\varphi_{n}(1-\varphi_{n})\frac{|u_n|^2}{|x|^2}dx+o(1)\\
&\leq& S_{\lambda,s}(\RT)+r_n^{-2}\int_{\RT}u_n^2dx+o(1)\\
&=& S_{\lambda,s}(\RT)+o(1)
\end{eqnarray*}
as $n\to\infty$. Hence $\alpha^{\frac{1}{3-s}}+(1-\alpha)^{\frac{1}{3-s}}\leq1$, which implies that $\alpha=1$ since $\frac12\leq\alpha\leq 1$. This proves the claim in \eqref{2.5}.

We now claim that there exist $u_\infty\in \X$ satisfying $u_n\rightharpoonup u_\infty$ in $\X$ as $n\to\infty$ and $x_0\not=0$ such that
\begin{eqnarray}
&\displaystyle {\rm either}\ \lim_{n\to\infty} \frac{|u_n|^{6-2s}}{|x|^s}dx=\frac{|u_\infty|^{6-2s}}{|x|^s}dx\ {\rm and}\ \int_{\RN}\frac{|u_\infty|^{6-2s}}{|x|^s}dx=1,\label{2.8}\\
&\displaystyle  {\rm or}\ \lim_{n\to\infty} \frac{|u_n|^{6-2s}}{|x|^s}dx=\delta_{x_0}\ {\rm and}\ u_\infty=0.\label{2.9}
\end{eqnarray}
Arguing as above, we get that for  all $x\in \RN$, we have that
\begin{eqnarray*}
\lim_{r\to0^+}\lim_{n\to\infty}\int_{B_r(x')\times B_r(x_3)}\frac{|u_n|^{6-2s}}{|x|^s} dx=\alpha_x\in\{0,1\}.
\end{eqnarray*}
Then follows from the second identity of \eqref{2.4} that $\alpha_0\leq \frac12$ and therefore $\alpha_0=0$. Moreover, it follows from the first identity of \eqref{2.4} that there exists at most one point $x_0\in\RT$ such that $\alpha_{x_0}=1$. In particular $x_0\not=0$ since  $\alpha_0=0$. It then follows from Lions's second concentration compactness lemma  that, up to a subsequence, there exist $u_\infty\in \X, x_0\in \RT\setminus\{0\}$ and $\nu\in\{0,1\}$ such that $u_n\rightharpoonup u_\infty$ weakly in $\X$ and
\begin{eqnarray*}
\lim_{n\to\infty} \frac{|u_n|^{6-2s}}{|x|^s}dx=\frac{|u_\infty|^{6-2s}}{|x|^s}dx+\nu\delta_{x_0}\ {\rm in\ the\ sense\ of\ measures}.
\end{eqnarray*}
In particular, due to \eqref{2.4} and \eqref{2.5}, we have that
\begin{eqnarray*}
1=\lim_{n\to\infty} \int_{\RT}\frac{|u_n|^{6-2s}}{|x|^s}dx= \int_{\RT} \frac{|u_\infty|^{6-2s}}{|x|^s}dx+\nu.
\end{eqnarray*}
Since $\nu\in \{0,1\}$, the claims in \eqref{2.8} and \eqref{2.9} follow.

We now assume that $u_\infty\not=0$, and we claim that $\lim_{n\to\infty}u_n=u_\infty$ strongly in $\X$ and that $u_\infty$ is an extremal for $S_{\lambda,s}(\RT)$.

Indeed, it follows from \eqref{2.8} that $\int_{\RT} \frac{|u_\infty|^{6-2s}}{|x|^s}dx=1$, hence
\begin{eqnarray}
S_{\lambda,s}(\RT)\leq\int_{\RT}\left(|\nabla u_\infty|^2+\frac{|u_\infty|^2}{|x'|^2}-\lambda\frac{|u_\infty|^2}{|x|^2}\right)dx.\label{2.10}
\end{eqnarray}
Moreover, since $u_n\rightharpoonup u_\infty$ weakly as $n\to\infty$, we have that
\begin{eqnarray}
\left.
  \begin{array}{ll}
&\displaystyle \int_{\RT}\left(|\nabla u_\infty|^2+\frac{|u_\infty|^2}{|x'|^2}-\lambda\frac{|u_\infty|^2}{|x|^2}\right)dx\\
&\displaystyle  \leq  \liminf_{n\to\infty} \int_{\RT}\left(|\nabla u_n|^2+ \frac{|u_n|^2}{|x'|^2}-\lambda\frac{|u_n|^2}{|x|^2}\right)dx\\
&\displaystyle  =S_{\lambda,s}(\RT).
  \end{array}
\right.\label{2.11}
\end{eqnarray}
Hence, combining with \eqref{2.10} and \eqref{2.11}, $u_\infty$ is an extremal for $S_{\lambda,s}(\RT)$ and boundedness yields the weak convergence of $u_n$ to $u_\infty$ in $\X$, furthermore, the fact $\lim_{n\to\infty}u_n=u_\infty$ strongly in $\X$ holds. This proved the claims.

We now assume $u_\infty\equiv0$. According to the fact $u_n\rightharpoonup u_\infty\equiv0$ weakly in $\X$ as $n\to\infty$, then for any $1\leq q<2^*(0), u_n\to0$ strongly in $L^q_{\rm loc}(\RT)$ when $n\to\infty$. It follows from $s>0$ that $2^*(s)<2^*(0)$, for $x_0\not=0$, we have that
\begin{eqnarray*}
\lim_{n\to\infty}\int_{B_\delta(x_0')\times B_\delta(x_{03})}\frac{|u_n|^{6-2s}}{|x|^s}dx=0
\end{eqnarray*}
for $\delta>0$ small enough, contradicting \eqref{2.9}. As a result that $u_\infty\not=0$.

Based on the proof above, it can be concluded that there exists a $u_{\lambda,s}\in \X$ such that
\begin{eqnarray*}
S_{\lambda,s}(\RT)=\frac{\int_{\RT}|\nabla u_{\lambda,s}|^2+\frac{|u_{\lambda,s}|^2}{|x'|^2}-\lambda \frac{|u_{\lambda,s}|^2}{|x|^2}dx}{\left(\int_{\RT}\frac{|u_{\lambda,s}|^{6-2s}}{|x|^s}\right)^{\frac{1}{3-s}}}.
\end{eqnarray*}
According to \eqref{1.14}, there exists a $\hu_{\lambda,s}\in \DF$ define by \eqref{1.16} replacing $u$ with $u_{\lambda,s}$ above such that
\begin{eqnarray*}
\frac{\int_{\RT}|\nabla\times \hu_{\lambda,s}|^2-\lambda \frac{|\hu_{\lambda,s}|^2}{|x|^2}dx}{\left(\int_{\RT}\frac{|\hu_{\lambda,s}|^{6-2s}}{|x|^s}\right)^{\frac{1}{3-s}}}=S_{\lambda,s}(\RT)=\inf_{\hu\in \DF}\frac{\int_{\RT}|\nabla \times\hu|^2-\lambda \frac{|\hu|^2}{|x|^2}dx}{\left(\int_{\RT}\frac{|\hu|^{6-2s}}{|x|^s}\right)^{\frac{1}{3-s}}}.
\end{eqnarray*}
Making a  scaling for $\hu_{\lambda,s}$, we   get a nontrivial ground state solution of \eqref{E}. \hbx

\section{The case  $\lambda\in \mathbb{R}\setminus\{0\}, 0\leq s_1<s_2<2$}\label{sec3}

In this section, we main consider the case $0\leq s_1<s_2<2, \lambda\in \mathbb{R}\setminus\{0\}$ and prove the Theorems \ref{thm1.2} and \ref{thm1.3} based on the Nehari manifold and the results for the case of $\lambda=0$ in Theorem \ref{thm1.1}. Let us simply denote $S_{0,s}(\RT)$ as $S_{s}$.  We first establish the key lemma which is important to show that the least energy is equal to the level of the mountain pass.
\begin{lem}\label{lem3.1} Assume $0\leq s_1<s_2<2, \lambda\in \mathbb{R}\setminus\{0\}$ hold.  For each $u\in \X\setminus\{0\}$, there exists a unique $t_u>0$ such that $t_u u\in \mathcal{N}$ and $\mathcal{I}_\lambda(t_u u)=\max_{t \geqslant 0} \I_\lambda(tu)$.  The function $u \mapsto t_u$ is continuous and the map $u \mapsto t_u u$ is a homeomorphism of the unit sphere in $\X$ with $\mathcal{N}$.
\end{lem}
\pf \  For any $u\in \X\setminus\{0\}, t>0$,
\begin{eqnarray*}
\frac{d \I_\lambda(tu)}{dt}=t\mathcal{A}(u)-t^{5-2s}\mathcal{B}(u)-t^{5-2s}\lambda\mathcal{C}(u).
\end{eqnarray*}
The fact $0\leq s_1<s_2<2$ implies that there exists a unique $t_u>0$ such that $\frac{d \I_\lambda(tu)}{dt}|_{t=t_u}=0$, that is $t_u u\in \mathcal{N}$. \cite[Chapter 4]{1996Willem} can be referenced for the remaining proof.\hfill $\Box$

We define\begin{eqnarray}
&& \bar{c}_\lambda:=\inf_{u\in  \X\setminus\{0\}}\max_{t\geqslant0}\I_\lambda(tu),\label{3.1}
\\
&&\hat{c}_\lambda:=\inf_{\gamma\in  \Gamma}\max_{t\in [0,1]}\I_\lambda(\gamma(t)),\label{3.2}
\end{eqnarray}
where
\begin{eqnarray}
\Gamma:=\left\{\gamma\in C([0,1],\X): \ \gamma(0)=0, \  \I_\lambda(\gamma(1))<0\right\}.\label{3.3}
\end{eqnarray}
By \cite[Chapter 4]{1996Willem} and Lemma \ref{lem3.1}  we have that
 \begin{eqnarray}
 m_\lambda=\bar{c}_\lambda=\hat{c}_\lambda.\label{3.4}
 \end{eqnarray}
As $\lambda=0$, the ground state solutions have been in Section \ref{sec3}. We will  prove the existence of ground state solutions in two different cases as $\lambda\not=0$.

Now we introduce a interpolation inequality which is a changed version of Caffarelli-Kohn-Nirenberg inequality after a suitable transform of functions in \cite{1984CKN}.
\begin{lem}\label{lemCKN}
Assume $0\leq s_1,s_2<2$. There exists a constant $\bar{S}$ such that for any $u\in \X$,
 \begin{eqnarray*}
\left(\int_{\mathbb{R}^3}\frac{|u|^{6-2s_1}}{|x|^{s_1}}dx\right)^{\frac{1}{6-2s_1}}\leq \bar{S}\left(\int_{\mathbb{R}^3}|\nabla u|^2+\frac{|u|^2}{|x'|^2}dx\right)^{\frac{a}{2}}\left(\int_{\mathbb{R}^3} \frac{|u|^{6-2s_2}}{|x|^{s_2}}dx\right)^{\frac{1-a}{6-2s_2}},
\end{eqnarray*}
where
\begin{eqnarray*}
1> a\geq \left\{
          \begin{array}{ll}
           \frac{3(s_2-s_1)}{s_2(3-s_1)} \ &{\rm if}\ 2\geq s_2>s_1>0,\\[1em]
  \frac{s_1-s_2}{(2-s_2)(3-s_1)} \  &{\rm if}\ 2>s_1>s_2\geq0.
          \end{array}
        \right.
\end{eqnarray*}
\end{lem}

\subsection{The case $\lambda>0, 0\leq s_1<s_2<2$}\label{sec3.1}
In this section, we consider the case of $\lambda>0$ and we always assume $0\leq s_1<s_2<2$.

 \begin{lem}\label{lem3.3}
There holds
\begin{eqnarray*}
0<\hat{c}_\lambda<c^*_\lambda:=\min\left\{\frac{2-s_1}{6-2s_1} S_{s_1}^{\frac{N-s_1}{2-s_1}},
\frac{2-s_2}{6-2s_2}\lambda^{-\frac{1}{2-s_2}} S_{s_2}^{\frac{3-s_2}{2-s_2}}\right\}
\end{eqnarray*} and there exists a ${\rm (PS)_{\hat{c}_\lambda}}$ sequence
$\{u_n\}\in \X$ of $\I_\lambda$.
\end{lem}
\pf
 It is obvious that $\I_\lambda(0)=0$.  On one hand, by the inequality \eqref{Inequality}, for any $u\in \X\setminus\{0\}$, it holds that
\begin{eqnarray*}
\I_\lambda(u)\geqslant\frac{1}{2}\|u\|^2-\frac{1}{6-2s_1}S_{s_1}^{s_1-3}
\|u\|^{6-2s_1}-\frac{\lambda}{6-2s_2}S_{s_2}^{s_2-3}
\|u\|^{6-2s_2}.
\end{eqnarray*}
It will be seen from this that there exists $\rho>0$ such that
\begin{eqnarray*}
\I_\lambda(u)\geqslant \varpi> 0 \ \ {\rm as} \ \ \|u\|=\rho,
\end{eqnarray*}
where
\begin{eqnarray*}
\varpi=\frac{1}{2}\rho^2-\frac{1}{6-2s_1}S_{s_1}^{s_1-3}
\rho^{6-2s_1}-\frac{\lambda}{6-2s_2}S_{s_2}^{s_2-3}
\rho^{6-2s_2}.
\end{eqnarray*}
On the other hand, for any fixed  $u\in \X\setminus\{0\}$,
\begin{eqnarray*}
\I_\lambda(tu)=\frac{t^2}{2}\mathcal{A}(u)-\frac{t^{6-2s_1}}{6-2s_1}\mathcal{B}(u)
-\frac{\lambda t^{6-2s_2}}{6-2s_2}\mathcal{C}(u) \to-\infty \ {\rm as} \ t \to \infty.
\end{eqnarray*}
 Thus there exists  $\upsilon\in\X$ satisfying
\begin{eqnarray*}
 \|\upsilon\|>\rho, \ \ \ \I_\lambda(\upsilon)<0.
\end{eqnarray*}
Now applying the mountain pass theorem, we obtain a ${\rm(PS)}$ sequence $\{u_n\}\subset \X$ of $\I_\lambda$ at the level   $\hat{c}_\lambda$.

 Next we prove that
 \begin{eqnarray*}
\hat{c}_\lambda<c^*_\lambda.
\end{eqnarray*}
By \eqref{3.4}  we only need to prove that
\begin{eqnarray*}
\bar{c}_\lambda<c^*_\lambda.
\end{eqnarray*}
Choosing the extremum function $u_{0,s_1}$ in \eqref{Inequality} with $s=s_1$, then we get that
\begin{eqnarray}
\left.
  \begin{array}{ll}
\displaystyle\max_{t\geqslant0}\I_\lambda(tu_{0,s_1})&=\frac{t_{s_1}^2}{2}\mathcal{A}(u_{0,s_1})
-\frac{t_{s_1}^{6-2s_1}}{6-2s_1}\mathcal{B}(u_{0,s_1})
-\frac{\lambda t_{s_1}^{6-2s_2}}{6-2s_2}\mathcal{C}(u_{0,s_1})\\
&\leqslant\displaystyle \max_{t\geqslant0}\left\{\frac{t^2}{2}\mathcal{A}(u_{0,s_1})
-\frac{t^{6-2s_1}}{6-2s_1}\mathcal{B}(u_{0,s_1})\right\}
-\frac{\lambda t_{s_1}^{6-2s_2}}{6-2s_2}\mathcal{C}(u_{0,s_1})\\
&=\displaystyle  \frac{2-s_1}{6-2s_1}S_{\alpha_1}^{\frac{3-s_1}{2-s_1}}
-\frac{\lambda t_{s_1}^{6-2s_2}}{6-2s_2}\mathcal{C}(u_{0,s_1})\\
& \displaystyle  <\frac{2-s_1}{6-2s_1}S_{\alpha_1}^{\frac{3-s_1}{2-s_1}},
  \end{array}
\right.\label{3.5}
\end{eqnarray}
where we have used a  fact that
\begin{eqnarray*}
\mathcal{C}(u_{0,s_1})=\int_{\mathbb{R}^N}\frac{|u_{0,s_1}|^{6-2s_2}}{|x|^{s_2}}dx>0.
\end{eqnarray*}
Similarly, taking the extremum function $u_{0,s_2}$ in \eqref{Inequality} with $s=s_2$, then we get that
\begin{eqnarray}
\displaystyle \max_{t\geqslant0}\I_\lambda(tu_{0,s_2})<\frac{2-s_2}{6-2s_2}\lambda^{-\frac{1}{2-s_2}}
S_{s_2}^{\frac{3-s_2}{2-s_2}}.\label{3.6}
\end{eqnarray}
Combining with \eqref{3.5} and \eqref{3.6}, we get that $\hat{c}_\lambda=\bar{c}_\lambda<c^*_\lambda$.   \hfill$\Box$

In the following lemmas we investigate the properties of the $\mathrm{(PS)}_{\hat{c}_\lambda}$ sequence $\{u_n\}$ of $\I_\lambda$ found in Lemma \ref{lem3.3}.

 \begin{lem}\label{lem3.4}
 If $u_n\rightharpoonup 0$ in $\X$, then for any domain $B_{a,b}(0)$, up to a subsequence and still denoted by $\{u_n\}$ such that
\begin{eqnarray}
\left.
  \begin{array}{ll} \displaystyle
\int_{B_{a,b}(0)}|\nabla u_n|^2+\frac{|u_n|^2}{|x'|^{2}}dx\to 0, \ \ \int_{B_{a,b}(0)}\frac{|u_n|^{6-2s_i}}{|x|^{s_i}}dx\to 0, \ \ \ i=1,2.
  \end{array}
\right.\label{3.7}
\end{eqnarray}
\end{lem}
{\bf Proof.}\ For any $R>r>0$, the compactness of the embedding
\begin{eqnarray}
\X\hookrightarrow\hookrightarrow L^{6-2s_2}(B_{r,R}(0);|x|^{-s_2})\label{3.8}
\end{eqnarray}
implies that
\begin{eqnarray*}
\int_{B_{a,b}(0)}\frac{|u_n|^{6-2s_2}}{|x|^{s_2}}dx \to 0.
\end{eqnarray*}
Let $\eta=\eta_1\eta_2$, where $\eta_1\in C_{0,r}^\infty(\mathbb{R}^2)$ such that $0\leqslant \eta_1\leqslant1, \eta_1(0')=0$ and $\eta_1|_{B_{a,b}(0')}\equiv1$, $\eta_2\in C_{0,r}^\infty(\mathbb{R})$ such that $0\leqslant \eta_2\leqslant1, \eta_2(0_3)=0$ and $\eta_2|_{B_{a,b}(0_3)}\equiv1$. Since $\eta^2u_n\in \X$ for all $n\in \mathbb{N}$, combining with \eqref{3.8}, we get that
\begin{eqnarray}
\left.
  \begin{array}{ll}
o(1)&=\langle\I_\lambda'(u_n),\eta^2u_n\rangle\\
&\displaystyle =\int_{\mathbb{R}^3}\nabla u_n\cdot \nabla(\eta^2u_n)+\frac{|\eta u_n|^2}{|x'|^{2}}dx
-\int_{\mathbb{R}^3}\frac{|u_n|^{6-2s_1}\eta^2}{|x|^{s_1}}dx-\lambda\int_{\mathbb{R}^3}\frac{|u_n|^{6-2s_2}\eta^2}{|x|^{s_2}}dx\\[1em]
&= \displaystyle \int_{\mathbb{R}^3}\eta^2|\nabla u_n|^2+2\eta u_n\nabla\eta\nabla u_n+\frac{|\eta u_n|^2}{|x'|^{2}}dx  -\int_{\mathbb{R}^3}\frac{|u_n|^{6-2s_1}\eta^2}{|x|^{s_1}}dx.
  \end{array}
\right.\label{3.9}
\end{eqnarray}
We claim that
\begin{eqnarray}
\int_{\mathbb{R}^3}2\eta u_n\nabla\eta\nabla u_ndx=o(1).\label{3.10}
\end{eqnarray}
Based the Sobolev embedding theorem, we can obtain that
\begin{eqnarray*}
\int_{{\rm supp}|\eta|}|u_n|^2dx\to0 \ {\rm as} \ n\to\infty,
\end{eqnarray*}
combining with H\"{o}lder inequality, then
\begin{eqnarray*}
\left|\int_{\mathbb{R}^3}2\eta u_n\nabla\eta\cdot\nabla u_ndx\right|\leq 2\left(\int_{\mathbb{R}^3}|\nabla\eta\cdot\nabla u_n|^2dx\right)^{\frac12}\left(\int_{\mathbb{R}^3}|\eta u_n|^2 dx\right)^{\frac12}\to 0 \ {\rm as} \ n\to\infty.
\end{eqnarray*}
Thus, by \eqref{3.9} and \eqref{3.10}, it is easy to see that
\begin{eqnarray*}
\int_{\mathbb{R}^3}|\eta\nabla u_n|^2+\frac{|\eta u_n|^2}{|x'|^{2}}dx-\int_{\mathbb{R}^3}\frac{|u_n|^{6-2s_1}\eta^2}{|x|^{s_1}}dx\to0\ {\rm as}\ n\to\infty,
\end{eqnarray*}
which implies that
\begin{eqnarray*}
\int_{\mathbb{R}^3}|\eta\nabla u_n|^2+\frac{|\eta u_n|^2}{|x'|^{2}}dx=\int_{\mathbb{R}^3}\frac{|u_n|^{6-2s_1}\eta^2}{|x|^{s_1}}dx+o(1)\ {\rm as}\ n\to\infty.
\end{eqnarray*}
Using the H\"{o}lder inequality, \eqref{Inequality} and \eqref{3.9}, \eqref{3.10}, we have, as $n\to\infty$,
\begin{eqnarray*}
\int_{\mathbb{R}^3}|\eta\nabla u_n|^2+\frac{|\eta u_n|^2}{|x'|^{2}}dx
&\leq&\left(\int_{\mathbb{R}^3}\frac{|u_n|^{6-2s_1}}{|x|^{s_1}}dx\right)^{\frac{2-s_1}{3-s_1}}
\left(\int_{\mathbb{R}^3}\frac{|\eta u_n|^{6-2s_1}}{|x|^{s_1}}dx\right)^{\frac{1}{3-s_1}}+o(1)\\
&\leq& \left(\int_{\mathbb{R}^3}\frac{|u_n|^{6-2s_1}}{|x|^{s_1}}dx\right)^{\frac{2-s_1}{3-s_1}}
S_{s_1}^{-1}\int_{\mathbb{R}^3}|\eta\nabla u_n|^2+\frac{|\eta u_n|^2}{|x'|^{2}}dx+o(1),
\end{eqnarray*}
it follows that
\begin{eqnarray}\label{3.11}
\left[1-\left(\int_{\mathbb{R}^N}\frac{|u_n|^{6-2s_1}}{|x|^{s_1}}dx\right)^{\frac{2-s_1}{3-s_1}}
S_{s_1}^{-1}\right]\int_{\mathbb{R}^N}|\eta\nabla u_n|^2+\frac{|\eta u_n|^2}{|x'|^{2}}dx\leq o(1)\ {\rm as}\ n\to\infty.
\end{eqnarray}
Since $\{u_n\}$ is a ${\rm (PS)_{\hat{c}_\lambda}}$ sequence, it is easy to see that, as $n\to\infty$,
\begin{eqnarray*}
\hat{c}_\lambda+o(1)&=&\I_\lambda(u_n)-\frac{1}{2}\langle\I_\lambda'(u_n), u_n\rangle\\
&=&\frac{2-s_1}{6-2s_1}\int_{\mathbb{R}^3}\frac{|u_n|^{6-2s_1}}{|x|^{s_1}}dx
+\lambda\frac{2-s_2}{6-2s_2}\int_{\mathbb{R}^3}\frac{|u_n|^{6-2s_2}}{|x|^{s_2}}dx,
\end{eqnarray*}
we can deduce that
\begin{eqnarray}\label{3.12}
\int_{\mathbb{R}^3}\frac{|u_n|^{6-2s_1}}{|x|^{s_1}}dx\leq \frac{6-2s_1}{4-s_1}\hat{c}_\lambda.
\end{eqnarray}
Combining with \eqref{3.11}, \eqref{3.12} and  $\hat{c}_\lambda<c^*_\lambda$  in Lemma \ref{lem3.3}, we have
\begin{eqnarray*}
\int_{\mathbb{R}^3}|\eta\nabla u_n|^2+\frac{|\eta u_n|^2}{|x'|^2}dx\to0\ {\rm as}\ n\to\infty.
\end{eqnarray*}
Since $\eta|_{B_{a,b}(0)}\equiv1$,  we get
\begin{eqnarray*}
\int_{B_{a,b}(0)}|\nabla u_n|^2+\frac{|\eta u_n|^2}{|x'|^2}dx\to0\ {\rm as}\ n\to\infty.
\end{eqnarray*}
We complete the proof. \hbx

For any $\delta>0$, we define
\begin{eqnarray}
\left.
  \begin{array}{ll}
 &\displaystyle \kappa_1:=\limsup_{n\to\infty}\int_{B_\delta(0')\times B_\delta(0_3)}\frac{|u_n|^{6-2s_1}}{|x|^{s_1}}dx,\\
& \displaystyle \kappa_2:=\limsup_{n\to\infty}\int_{B_\delta(0')\times B_\delta(0_3)}\frac{|u_n|^{6-2s_2}}{|x|^{s_2}}dx,\\
&\displaystyle \kappa :=\limsup_{n\to\infty}\int_{B_\delta(0')\times B_\delta(0_3)}|\nabla u_n|^2+\frac{|u_n|^2}{|x'|^{2}}dx.
  \end{array}
\right.\label{3.13}
\end{eqnarray}
It follows from Lemma \ref{lem3.4} that these three quantities are well defined and independent of the choice of $\delta>0$.

\begin{lem}\label{lem3.5}
If $u_n\rightharpoonup 0$ in $\X$, then there exist $\epsilon_0:=\epsilon_0(s_1,s_2,\hat{c}_\lambda,\lambda)>0$ and subsequence{\rm(}still denoted by $\{u_n\}${\rm)} such that
\begin{eqnarray*}
{\rm either} \ \ \lim_{n\to\infty}\int_{B_\delta(0')\times B_\delta(0_3)}\frac{|u_n|^{6-2s_1}}{|x|^{s_1}}dx=0 \ \ \ {\rm or} \ \
\lim_{n\to\infty}\int_{B_\delta(0')\times B_\delta(0_3)}\frac{|u_n|^{6-2s_1}}{|x|^{s_1}}dx\geqslant\epsilon_0
\end{eqnarray*}
for all $\delta>0$.
\end{lem}
 \pf Let $\phi=\phi_1\phi_2$, where $\phi_1\in C_{0,r}^\infty(\mathbb{R}^2)$ is nonnegative and satisfy $\phi_1|_{B_\delta(0')}\equiv1$ with $\delta>0$, and $\phi_2\in C_{0,r}^\infty(\mathbb{R})$ is nonnegative and  satisfy $\phi_2|_{B_\delta(0_3)}\equiv1$ with $\delta>0$. It follows from $\phi u_n\in \X$ that as $n\to\infty$,
\begin{eqnarray}\label{3.14}
\int_{\mathbb{R}^3}\nabla u_n\nabla (\phi u_n)+\frac{\phi|u_n|^2}{|x'|^{2}}dx
-\int_{\mathbb{R}^3}\frac{|u_n|^{6-2s_1}\phi }{|x|^{s_1}}dx-\lambda\int_{\mathbb{R}^3}\frac{|u_n|^{6-2s_2}\phi}{|x|^{s_2}}dx\to0.
\end{eqnarray}
By \eqref{3.7} in Lemma \ref{lem3.4}, we obtain that
\begin{eqnarray*}
 \int_{\mathbb{R}^3}\nabla u_n\nabla (\phi u_n)+\frac{\phi|u_n|^2}{|x'|^{2}}dx
&=&\int_{\mathbb{R}^3}\phi |\nabla u_n|^2+u_n\nabla u_n\nabla\phi+\frac{\phi|u_n|^2}{|x'|^{2}}dx\\
&=&\int_{B_\delta(0')\times B_\delta(0_3)}|\nabla u_n|^2+\frac{|u_n|^2}{|x'|^{2}}dx+o(1),
\end{eqnarray*}
\begin{eqnarray*}
\int_{\mathbb{R}^3}\frac{|u_n|^{6-2s_i}\phi}{|x|^{s_i}}dx
\to \int_{B_\delta(0')\times B_\delta(0_3)}\frac{|u_n|^{6-2s_i}}{|x|^{s_i}}dx\ {\rm as}\ n\to\infty, i=1,2.
\end{eqnarray*}
The limit \eqref{3.14} implies that
\begin{eqnarray}
\kappa=\kappa_1+\lambda\kappa_2.\label{3.15}
\end{eqnarray}
The definition \eqref{Inequality} leads to
\begin{eqnarray*}
\left(\int_{\mathbb{R}^3}\frac{|\phi u_n|^{6-2s_1}}{ |x|^{s_1}}dx\right)^{\frac{1}{3-s_1}}\leqslant S_{s_1}^{-1}\int_{\mathbb{R}^3}|\nabla (\phi u_n)|^2+\frac{\phi|u_n|^2}{|x'|^2}dx.
\end{eqnarray*}
Thus
\begin{eqnarray*}
\left(\int_{B_\delta(0')\times B_\delta(0_3)} \frac{|u_n|^{6-2s_1}}{ |x|^{s_1}}dx\right)^{\frac{1}{3-s_1}}\leqslant S_{s_1}^{-1}\int_{B_\delta(0')\times B_\delta(0_3)}|\nabla u_n|^2+\frac{\phi|u_n|^2}{|x'|^{2}}dx\ {\rm as}\ n\to\infty.
\end{eqnarray*}
Furthermore,
\begin{eqnarray}
\kappa_1^{\frac{1}{3-s_1}}\leqslant S_{s_1}^{-1}\kappa.\label{3.16}
\end{eqnarray}
The conclusions \eqref{3.15} and \eqref{3.16} lead to
\begin{eqnarray*}
\kappa_1^{\frac{1}{3-s_1}}\leqslant S_{s_1}^{-1}\kappa= S_{s_1}^{-1}\kappa_1+S_{s_1}^{-1}\lambda\kappa_2.
\end{eqnarray*}
It follows that
\begin{eqnarray}
\kappa_1^{\frac{1}{3-s_1}}\left(1-S_{s_1}^{-1}
\kappa_1^{\frac{2-s_1}{3-s_1}}\right)\leqslant S_{s_1}^{-1}\lambda\kappa_2.\label{3.17}
\end{eqnarray}
Since $\{u_n\}$ is a bounded $(\mathrm{PS})_{\hat{c}_\lambda}$ sequence, it is obvious that
\begin{eqnarray*}
\I_\lambda(u_n)-\frac{1}{2}\langle\I_\lambda'(u_n),u_n\rangle=\frac{2-s_1}{6-2s_1}
\int_{\mathbb{R}^3}\frac{|u_n|^{6-2s_1}}{|x|^{s_1}}dx+\frac{\lambda(2-s_2)}{6-2s_2}
\int_{\mathbb{R}^3}\frac{|u_n|^{6-2s_2}}{|x|^{s_2}}dx=\hat{c}_\lambda+o(1).
\end{eqnarray*}
It is easy to see that
\begin{eqnarray*}
\int_{\mathbb{R}^3}\frac{|u_n|^{6-2s_1}}{|x|^{s_1}}dx
\leqslant\frac{6-2s_1}{2-s_1} \hat{c}_\lambda+o(1)
\end{eqnarray*}
and
\begin{eqnarray}
\kappa_1\leqslant\frac{6-2s_1}{2-s_1}\hat{c}_\lambda.\label{3.18}
\end{eqnarray}
Combining with \eqref{3.17} and \eqref{3.18}, we deduce that
\begin{eqnarray*}
\kappa_1^{\frac{1}{3-s_1}}\left(1-S_{s_1}^{-1}\left(\frac{6-2s_1}{2-s_1}\hat{c}_\lambda \right)^{\frac{2-s_1}{3-s_1}}\right)\leqslant S_{s_1}^{-1}\lambda\kappa_2.
\end{eqnarray*}
The fact $\hat{c}_\lambda<c^*_\lambda$ implies that
\begin{eqnarray*}
S_{s_1}^{-1}\left(\frac{6-2s_1}{2-s_1}\hat{c}_\lambda\right)^
{\frac{2-s_1}{3-s_1}}<1.
\end{eqnarray*}
Thus there exists $\delta_1>0$ depending on $s_1,S_{s_1}, \hat{c}_\lambda$ such that $\kappa_1^{\frac{1}{3-s_1}}\leqslant \delta_1\kappa_2$. Similarly, we have $\kappa_2^{\frac{1}{3-s_2}}\leqslant \delta_2\kappa_1$ for some $\delta_2>0$.  It follows that  there exists $\epsilon_0:=\epsilon_0(s_1,s_2, \hat{c}_\lambda,\lambda)>0$ such that
\begin{eqnarray*}
{\rm either} \ \ \kappa_1=\kappa_2=0 \ \ \ {\rm or} \ \ \ \kappa_1\geqslant \epsilon_0, \ \ \kappa_2\geqslant\epsilon_0.
\end{eqnarray*}
We complete the proof. \hfill $\Box$

We define
 \begin{eqnarray}
 \tilde{u}_n(x): = r_n^{\frac{1}{2}}u_n (r_nx)  \ \  {\rm for} \ \ x\in \mathbb{R}^3,   \ r_n>0.\label{3.19}
\end{eqnarray}
Then $\{\tilde{u}_n\}\subset \X$ is also a ${\rm(PS)_{\hat{c}_\lambda}}$ sequence of $\I_\lambda$.

\begin{lem}\label{lem3.6}  There exists $\epsilon_1\in (0,\frac{\epsilon_0}{2}]$ such that for all $\epsilon\in(0, \epsilon_1)$, there exists  a sequence $\{r_n>0\}$
 such that  $\{\tilde{u}_n\}$  verifies
\begin{eqnarray}
\int_{B_1(0')\times B_1(0_3)}\frac{|\tilde{u}_n|^{6-2s_1}}{|x|^{s_1}}dx=\epsilon. \label{3.20}
\end{eqnarray}
\end{lem}
\pf  Since $\hat{c}_\lambda>0$, it follows from the  inequality in Lemma \ref{lemCKN} that
 $$\mathcal{B}_\infty:=\lim_{n\to\infty}\int_{\mathbb{R}^3}\frac{|u_n|^{6-2s_1}}{|x|^{s_1}}dx>0.$$
Let $\epsilon_1:=\min\{\frac{\epsilon_0}{2},\mathcal{B}_\infty\}$. For fixed $\epsilon\in (0,\epsilon_1)$, up to a subsequence, still denoted by $\{u_n\}$,
 for any $n\in \mathbb{N}$,  there exists $r_n>0$ such that
\begin{eqnarray*}
\int_{B_{r_n}(0')\times B_{r_n}(0_3)}\frac{|u_n|^{6-2s_1}}{|x|^{s_1}}dx=\epsilon.
\end{eqnarray*}
By dilating transformation, it is easy to check that $\{\tilde{u}_n\}$ satisfies \eqref{3.20}.
 \hfill$\Box$

Now we are ready to give the proof of Theorem \ref{thm1.2} for the case of $\lambda>0, 0\leq s_1<s_2<2$.

\noindent{\bf Proof of Theorem \ref{thm1.2}.} \ By Lemma \ref{lem3.3},  $\I_\lambda$ has a  ${\rm(PS)_{\hat{c}_\lambda}}$ sequence $\{u_n\}\subset \X$. By Lemma \ref{lem3.6},  the sequence $\{\tilde{u}_n\}$ defined by \eqref{3.19}  satisfies \eqref{3.20} and is also  a  ${\rm(PS)_{\hat{c}_\lambda}}$ sequence of $\I_\lambda$.
 Thus
\begin{eqnarray}
\left.
  \begin{array}{ll} \displaystyle
 &\displaystyle   \I_\lambda(\tilde{u}_n)-\frac{1}{6-2s_2}\langle\I_\lambda'(\tilde{u}_n),\tilde{u}_n\rangle\\
& \displaystyle  \geqslant \frac{2-s_2}{6-2s_2}\int_{\mathbb{R}^3}|\nabla \tilde{u}_n|^2+\frac{|\tilde{u}_n|^2}{|x'|^{2}}dx+\frac{s_2-s_1}{2(3-s_1)(3-s_2)} \mathcal{B}(\tilde{u}_n)\\[3mm]
 &\displaystyle  \geqslant\frac{2-s_2}{6-2s_2}\int_{\mathbb{R}^3}|\nabla \tilde{u}_n|^2+\frac{|\tilde{u}_n|^2}{|x'|^{2}}dx.
  \end{array}
\right.\label{3.21}
\end{eqnarray}
It follows that $\{\tilde{u}_n\}$ is bounded in $\X$. Thus there is $\tilde{u}_0\in \X$ such that
\begin{eqnarray*}
  \left\{\begin{array}{ll}
\tilde{u}_n\rightharpoonup \tilde{u}_0 \ \ &\mbox{in} \ \X;\\
\tilde{u}_n\rightharpoonup \tilde{u}_0 \ \ &\mbox{in} \ L^{6-2s_i}(\mathbb{R}^3;|x|^{-s_i}), \ i=1, 2; \\
\tilde{u}_n(x)\rightarrow \tilde{u}_0(x) \ \ & \mbox{a.e. \ on} \  \mathbb{R}^3.
  \end{array}
\right.
\end{eqnarray*}
From the above it can be concluded  that $\tilde{u}_0$ is a solution of \eqref{1.17}, furthermore, copying the calculation process of \eqref{3.21}, we get
that $\I_\lambda(\tilde{u}_0)\geqslant0$.  Let $v_n:=\tilde{u}_n-\tilde{u}_0$.
Then $\{v_n\}$ is bounded in $\X$. Define
\begin{eqnarray*}
\mathcal{A}(v_n) \to \mathcal{A}_\infty, \ \  \mathcal{B}(v_n)\to \mathcal{B}_\infty, \ \ \mathcal{C}(v_n)\to \mathcal{C}_\infty.
\end{eqnarray*}
Then by Brezis-Lieb Lemma \cite{1983B-Lieb}, we have
\begin{eqnarray}
 && \I_\lambda(v_n)\to \frac{1}{2}\mathcal{A}_\infty-\frac{1}{6-2s_1}\mathcal{B}_\infty
-\frac{\lambda}{6-2s_2}\mathcal{C}_\infty=\hat{c}_\lambda-\I_\lambda(\tilde{u}_0),\label{3.22} \\
&&\langle\I_\lambda'(v_n),v_n\rangle\to \mathcal{A}_\infty-\mathcal{B}_\infty-\lambda\mathcal{C}_\infty=0.\label{3.23}
\end{eqnarray}
 If $\mathcal{A}_\infty=0$, then $\I_\lambda(\tilde{u}_0)=\hat{c}_\lambda$ and $\tilde{u}_0$ is a ground state solution of \eqref{1.17}. Assume $\mathcal{A}_\infty>0$ and $\tilde{u}_0=0$. Then  Lemma \ref{lem3.5} implies that either  $\displaystyle\lim_{n\to\infty}\int_{B_1(0')\times B_1(0_3)}\frac{|\tilde{u}_n|^{6-2s_1}}{|x|^{s_1}}dx=0$ or $\displaystyle\lim_{n\to\infty}\int_{B_1(0')\times B_1(0_3)}\frac{|\tilde{u}_n|^{6-2s_1}}{|x|^{s_1}}dx \geqslant \epsilon_0$. By Lemma \ref{lem3.6}, this is a contradiction to \eqref{3.20} as $0<\epsilon<\frac{\epsilon_0}{2}$. It must be $\tilde{u}_0\not=0$ and $\tilde{u}_0$ is a nontrivial solution of \eqref{E}.  If $\I_\lambda(\tilde{u}_0)=\hat{c}_\lambda$  then we complete the proof by \eqref{3.4}.  Otherwise, combining with the key fact \eqref{3.4}, we deduce that
\begin{eqnarray}
\I_\lambda(\tilde{u}_0)>\hat{c}_\lambda.\label{3.24}
\end{eqnarray}
Since
  \begin{eqnarray*}
 \I_\lambda(v_n)-\frac{1}{6-2s_2}\langle\I_\lambda'(v_n),v_n\rangle\geqslant\frac{2-s_2}{6-2s_2}\mathcal{A}(v_n)
\geqslant0,
\end{eqnarray*}
 it follows from \eqref{3.22} and \eqref{3.23} that
 \begin{eqnarray*}
\I_\lambda(\tilde{u}_0)\leqslant\hat{c}_\lambda,
\end{eqnarray*}
which is a contradiction with \eqref{3.24}. Thus $\tilde{u}_0$ is a ground state  solution of \eqref{1.17}.  The result  of Theorem \ref{thm1.2} is proved according to  Remark \ref{rem}.   \hfill $\Box$

\subsection{The case $\lambda<0, 0<s_1<s_2<2$}\label{sec3.2}
In this subsection, we may consider the case of $\lambda<0$ and  $0\leq s_1<s_2<2$. Applying the mountain pass theorem in \cite{1973AR}, we have the following lemma.
\begin{lem}\label{lem3.7}
Let $\lambda<0, 0<s_1<s_2<2$. There exists a sequence $\{u_n\}\subset \X$ such that
\begin{eqnarray}\label{3.25}
 \I_\lambda(u_n) \to \hat{c}_\lambda>0,  \ \ \ \ \ \I_\lambda'(u_n)\to0, \ \ \ n \to\infty
 \end{eqnarray}
with $\hat{c}_\lambda$ is in \eqref{3.2}.
\end{lem}

The properties of the $\mathrm{(PS)}_{\hat{c}}$ sequence $\{u_n\}$ of $\I_\lambda$ found in Lemma \ref{lem3.7}  will be investigated  in the following lemmas.

 \begin{lem}\label{lem3.8}
 If $u_n\rightharpoonup 0$ in $\X$, then for any  domain $B_{a,b}(0)$, up to a subsequence and still denoted by $\{u_n\}$ such that
\begin{eqnarray}
\left.
  \begin{array}{ll} \displaystyle
\int_{B_{a,b}(0)}|\nabla u_n|^2+\frac{|u_n|^2}{|x'|^{2}}dxdx\to 0, \ \ \int_{B_{a,b}(0)}\frac{|u_n|^{6-2s_i}}{|x|^{s_i}}dx\to 0, \ \ \ i=1,2.
  \end{array}
\right.\label{3.26}
\end{eqnarray}
\end{lem}
\pf  For any $R>r>0$, the compactness of the embedding
\begin{eqnarray}
\X\hookrightarrow\hookrightarrow L^{6-2s_i}(B_{r,R}(0);|x|^{-s_i}), \ i=1,2\label{3.27}
\end{eqnarray}
implies that
\begin{eqnarray*}
\int_{B_{a,b}(0)}\frac{|u_n|^{6-2s_i}}{|x|^{s_i}}dx \to 0  \ \ \ {\rm as}\ n\to\infty,\ i=1,2.
\end{eqnarray*}
Applying the function $\eta$ defined in Lemma \ref{lem3.4}.  Then $\eta^2u_n\in \X$ for all $n\in \mathbb{N}$, combining with \eqref{3.27} and \eqref{3.10}, we get that
\begin{eqnarray}
\left.
  \begin{array}{ll}
o(1)&=\langle\I_\lambda'(u_n),\eta^2u_n\rangle\\
&\displaystyle =\int_{\mathbb{R}^N}\nabla u_n\cdot \nabla(\eta^2u_n)+\frac{|\eta u_n|^2}{|x'|^{2}}dx
-\int_{\mathbb{R}^N}\frac{|u_n|^{6-2s_1}\eta^2}{|x|^{s_1}}dx
-\lambda\int_{\mathbb{R}^N}\frac{|u_n|^{6-2s_2}\eta^2}{|x|^{s_2}}dx\\
&= \displaystyle \int_{\mathbb{R}^N}\eta^2|\nabla u_n|^2+\frac{|\eta u_n|^2}{|x'|^{2}}dx
  \end{array}
\right.\label{445}
\end{eqnarray}
and
\begin{eqnarray*}
\int_{\mathbb{R}^N}|\eta\nabla u_n|^2+\frac{|\eta u_n|^2}{|x'|^{2}}dx\to0.
\end{eqnarray*}
Since $\eta|_{B_{a,b}(0)}\equiv1$, thus we get
\begin{eqnarray*}
\int_{B_{a,b}(0)}|\nabla u_n|^2+\frac{|u_n|^2}{|x'|^{2}}dx\to0\ {\rm as}\ n\to\infty.
\end{eqnarray*}
We complete the proof. \hbx

It follows from Lemma \ref{lem3.8} that those three quantities in \eqref{3.13} are well defined and independent of the choice of $\delta>0$.

\begin{lem}\label{lem3.9}
If $u_n\rightharpoonup 0$ in $\X$, then there exists  subsequence{\rm(}still denoted by $\{u_n\}${\rm)} such that
\begin{eqnarray*}
{\rm either} \ \ \lim_{n\to\infty}\int_{B_\delta(0')\times B_\delta(0_3)}\frac{|u_n|^{6-2s_1}}{|x|^{s_1}}dx=0 \ \ \ {\rm or} \ \
\lim_{n\to\infty}\int_{B_\delta(0')\times B_\delta(0_3)}\frac{|u_n|^{6-2s_1}}{|x|^{s_1}}dx\geqslant S_{s_1}^{\frac{3-s_1}{2-s_1}}
\end{eqnarray*}
for all $\delta>0$.
\end{lem}
 \pf Taking the function $\phi$ and coping the proof in Lemma \ref{lem3.5},  we obtain that
\begin{eqnarray}
\kappa=\kappa_1+\lambda\kappa_2.\label{3.29}
\end{eqnarray}
\begin{eqnarray}
\kappa_1^{\frac{1}{3-s_1}}\leqslant S_{s_1}^{-1}\kappa.\label{3.30}
\end{eqnarray}
The conclusions \eqref{3.29} and \eqref{3.30} with $\lambda<0$ lead to
\begin{eqnarray*}
\kappa_1^{\frac{1}{3-s_1}}\leqslant S_{s_1}^{-1}\kappa\leq S_{s_1}^{-1}\kappa_1.
\end{eqnarray*}
It follows that
\begin{eqnarray*}
\kappa_1=0\ \ {\rm or}\ \ \kappa_1\geq S_{s_1}^{\frac{3-s_1}{2-s_1}}.
\end{eqnarray*}
We complete the proof. \hfill $\Box$

Next we consider the transform \eqref{3.19}, then $\{\tilde{u}_n\}\subset \X$ is also a ${\rm(PS)_{\hat{c}_\lambda}}$ sequence of $\I_\lambda$.

\begin{lem}\label{lem3.10}  There exists $\epsilon_1\in (0,\frac{1}{2}S_{s_1}^{\frac{3-s_1}{2-s_1}}]$ such that for all $\epsilon\in(0, \epsilon_1)$, there exists  a sequence $\{r_n>0\}$
 such that  $\{\tilde{u}_n\}$  verifies
\begin{eqnarray}
\int_{B_1(0')\times B_1(0_3)}\frac{|\tilde{u}_n|^{6-2s_1}}{|x|^{s_1}}dx=\epsilon. \label{3.32}
\end{eqnarray}
\end{lem}
\pf  Since $\hat{c}_\lambda>0$, it follows from the interpolation inequality in Lemma \ref{lemCKN} that $$\mathcal{B}_\infty:=\lim_{n\to\infty}\int_{\mathbb{R}^3}\frac{|u_n|^{6-2s_1}}{|x|^{s_1}}dx>0.$$
Let $\epsilon_1\leq\frac{1}{2}S_{s_1}^{\frac{3-s_1}{2-s_1}}$. For fixed $\epsilon\in (0,\epsilon_1)$, up to a subsequence, still denoted by $\{u_n\}$,
 for any $n\in \mathbb{N}$,  there exists $r_n>0$ such that
\begin{eqnarray*}
\int_{B_{r_n}(0')\times B_{r_n}(0_3)}\frac{|u_n|^{6-2s_1}}{|x|^{s_1}}dx=\epsilon.
\end{eqnarray*}
By dilating transformation, it is easy to check that $\{\tilde{u}_n\}$ satisfies \eqref{3.32}.
 \hfill$\Box$

Based on the lemmas above, the results of Theorem \ref{thm1.3} with $\lambda<0, 0<s_1<s_2<2$ can be proved by coping the similar proof of Theorem \ref{thm1.2} in  subsection \ref{sec3.1}.

\section{The case $\lambda<0, 0\leq s_2<s_1<2$} \label{sec4}
In the present section, we focus on the proofs of Theorem \ref{thm1.4} and Theorem \ref{thm1.5}.

\subsection{Non-existence}\label{sec4.1}

In this subsection, we prove that the equation \eqref{E} has only zero solution in $\DF$ as $\lambda<\lambda^*$.

{\bf Proof of Theorem \ref{thm1.4}.}\  Based on the inequality in Lemma \ref{lemCKN} with $a=a_0:=\frac{s_1-s_2}{(2-s_2)(3-s_1)}$, we can directly obtain
\begin{eqnarray}
\left.
  \begin{array}{ll}
\displaystyle  \langle\I_\lambda'(u),u\rangle&\displaystyle =\int_{\mathbb{R}^3}|\nabla u|^2+\frac{|u|^2}{|x'|^{2}}dx-\int_{\mathbb{R}^3}\frac{|u|^{6-2s_1}}{|x|^{s_1}}dx-
\lambda\int_{\mathbb{R}^3}\frac{|u|^{6-2s_2}}{|x|^{s_2}}dx\\[1em]
\displaystyle &\geq \left(1-\bar{S}^{6-2s_1}a_0(3-s_1)\gamma^{\frac{1}{a_0(3-s_1)}}\right)\int_{\mathbb{R}^3}|\nabla u|^2+\frac{|u|^2}{|x'|^{2}}dx\\[1em]
&\displaystyle -\left(\lambda+\bar{S}^{6-2s_1}\frac{(1-a_0)(3-s_1)}{3-s_2}\gamma^{-\frac{(1-a_0)(3-s_1)}{(3-s_2)}}
\right)\int_{\mathbb{R}^3}\frac{|u|^{6-2s_2}}{|x|^{s_2}}dx
  \end{array}
\right.\label{4.1}
\end{eqnarray}
 and taking
\begin{eqnarray*}
\gamma=\frac{1}{2}\left[\left(\frac{(2-s_2)\bar{S}^{-(6-2s_1)}}{s_1-s_2}\right)^{\frac{s_1-s_2}{2-s_2}}
+\left(-\frac{2-s_1}{\lambda(2-s_2)}\bar{S}^{6-2s_1}\right)^{\frac{2-s_2}{2-s_1}}\right].
\end{eqnarray*}
Since $\lambda<\lambda^*$, we have
\begin{eqnarray*}
1-\bar{S}^{6-2s_1}a_0(3-s_1)\gamma^{\frac{1}{a_0(3-s_1)}}>0,\\ \lambda+\bar{S}^{6-2s_1}\frac{(1-a_0)(3-s_1)}{3-s_2}\gamma^{-\frac{(1-a_0)(3-s_1)}{3-s_2}}<0.
\end{eqnarray*}
Thus for any $u\in \X\setminus\{0\}$, as $\lambda<\lambda^*$, we have $\langle\I_\lambda'(u),u\rangle>0$, which implies the problem \eqref{1.17} has only zero solution. According to Remark \ref{rem}, we know that the problem \eqref{E} has only zero solution. And we complete the proof of Theorem \ref{thm1.4}. \hbx

\subsection{Existence}\label{sec4.2}
In this subsection, we prove that the equation \eqref{E} has solution as $\lambda$ small enough. According to Remark \ref{rem}, we need to prove the equation \eqref{1.17} has solution as $\lambda$ small enough. The assumptions $\lambda<0, 0\leq s_2<s_1<2$ make the undesirable  obstacle for establishing the mountain pass structure of equation \eqref{1.17}. The cut-off method in \cite{2006JS} is applied to overcome this obstacle. For any fixed $S>0$, the cut-off functional $\mathcal{J}_S: \X\to \mathbb{R}$ is defined as
 \begin{eqnarray}\label{4.2}
 \mathcal{J}_S(u)=\frac{1}{2}\int_{\mathbb{R}^3}|\nabla u|^2+\frac{|u|^2}{|x'|^2}dx
 -\frac{1}{6-2s_1}\int_{\mathbb{R}^3}\frac{|u|^{6-2s_1}}{|x|^{s_1}}dx
 -\frac{\lambda\Psi_S(u)}{6-2s_2} \int_{\mathbb{R}^3} \frac{|u|^{6-2s_2}}{|x|^{s_2}}dx,
 \end{eqnarray}
 where
  \begin{eqnarray}\label{4.3}
\Psi_S(u)=\psi\left(\frac{\|u\|^2}{S^2}\right)\geqslant0
  \end{eqnarray}
  and $\psi\in C_0^{\infty}(\mathbb{R},[0,1])$ satisfies $\psi(t)=1$ for $t\in[0,\frac{1}{2}]$ and supp $\psi\subset[0,1]$. The derivative of $\mathcal{J}_S$ is given by
 \begin{eqnarray*}
 \langle\mathcal{J}'_S(u),\varphi\rangle&=&\int_{\mathbb{R}^3} \nabla u\nabla \varphi+\frac{u\varphi}{|x'|^2} dx -\int_{\mathbb{R}^3}\frac{|u|^{4-2s_1}u\varphi}{|x|^{s_1}}dx-\lambda\Psi_S(u)\int_{\mathbb{R}^N}\frac{|u|^{4-2s_2}u\varphi}{|x|^{s_2}}dx\\
 &&\ \ -\frac{2\lambda\int_{\mathbb{R}^3} \nabla u\nabla\varphi+\frac{u\varphi}{|x'|^2} dx }{2^*(s_2)S^2}\psi'\left(\frac{\|u\|^2}{S^2}\right)\int_{\mathbb{R}^3} \frac{|u|^{6-2s_2}}{|x|^{s_2}}dx.
 \end{eqnarray*}
 We first show that the functional $\mathcal{J}_S$ has a mountain pass geometry for each fixed $S>0$.
\begin{lem} \label{lem4.1} Let  $\lambda<0$. We have
\begin{enumerate}
  \item[{\rm (i) }] there exist $\tilde{\rho}>0$ and $\tilde{\delta}>0$ such that $\mathcal{J}_S(u)\geqslant \tilde{\delta}$ for any $u\in \X$  with $\|u\|= \tilde{\rho}$;
     \item[{\rm (ii) }] there exists  $\tilde{v}\in \X$ satisfying $\|\tilde{v}\|> \tilde{\rho}$ and $\mathcal{I}_S(\tilde{v})<0.$
  \end{enumerate}
   \end{lem}
\pf (i) Since $\lambda<0$, by \eqref{4.3}, using the inequality \eqref{Inequality}, we derive that
\begin{eqnarray*}
 \mathcal{J}_S(u)\geqslant \frac{1}{2}\|u\|^2-S_{s_1}^{s_1-3}\|u\|^{6-2s_1}.
\end{eqnarray*}
It is not difficult to prove the conclusion (i) holds.

(ii) For any fixed $u\in \X\setminus\{0\}$ and $t>\frac{S}{\|u\|}$, we conclude  that
  \begin{eqnarray*}
  \mathcal{J}_S(tu)&=&\frac{t^2}{2}\int_{\mathbb{R}^3}|\nabla u|^2+\frac{|u|^2}{|x'|^2}dx
  -\frac{t^{6-2s_1}}{6-2s_1}\int_{\mathbb{R}^3}\frac{|u|^{6-2s_1}}{|x|^{s_1}}dx
 -\frac{\lambda t^{6-2s_2}}{6-2s_2}\Psi_S(tu) \int_{\mathbb{R}^3} \frac{|u|^{6-2s_2}}{|x|^{s_2}}dx\\
 &=&\frac{t^2}{2}\int_{\mathbb{R}^3}|\nabla u|^2+\frac{|u|^2}{|x'|^2}dx
 -\frac{t^{6-2s_1}}{6-2s_1}\int_{\mathbb{R}^3}\frac{|u|^{6-2s_1}}{|x|^{s_1}}dx.
 \end{eqnarray*}
 Taking $\tilde{v}=t_0 u$, where $t_0>\frac{S}{\|u\|}$ is large enough. Since $6-2s_1>2$, it is easy to see that (ii) holds.
 The proof is complete. $\hfill\Box$

  By Lemma \ref{lem4.1} and $\mathcal{J}_S(0)=0$, a mountain pass level for $\mathcal{J}_S$ can be defined as
  \begin{eqnarray}
  \bar c_S=\inf_{\gamma\in \Gamma} \max_{t\in [0, 1]} \mathcal{J}_S(\gamma(t))>0, \label{4.4}
  \end{eqnarray}
  where
  $$\Gamma_S =\left\{\gamma\in C([0, 1], \X): \ \gamma(0)=0, \ \gamma(1)=\tilde{v}\right\}.$$
Applying the mountain pass theorem,  there exists $\{u_n\} \subset \X$ satisfying
 \begin{eqnarray} \label{4.5}
\mathcal{J}_S(u_n) \rightarrow \bar{c}_S, \ \ \ \   \mathcal{J}_S'(u_n)\rightarrow 0, \ \ \mbox{as}\ n\rightarrow\infty.
\end{eqnarray}
Now we investigate the property of the sequence $\{u_n\}$ satisfying  \eqref{4.5}. We have
  \begin{lem}\label{lem4.2}  Let $\{u_n\} \subset \X$ satisfy \eqref{4.5}.
  Then for $S>0$ large enough, there exists $\lambda^{**}=\lambda^{**}(S)<0$ such that for any $\lambda<\lambda^{**}$,
 \begin{eqnarray}\label{4.6}
\limsup_{n\rightarrow\infty}\|u_n\|<\frac{S}{2}.
  \end{eqnarray}
 \end{lem}
 \pf  \ We first claim that $\{u_n\}$ is bounded. If $\|u_n\|\rightarrow\infty$ as
 $n\rightarrow\infty,$ then it follows from \eqref{4.3} that
  $$\Psi_S(u_n)=\psi\left(\frac{\|u_n\|^2}{S^2}\right)=0, \ \ \  \hbox{for  all large} \ \ n\in \mathbb{N}, $$
 and therefore for all $n\in \mathbb{N}$ large,
\begin{eqnarray*}
 \mathcal{J}_S(u_n)=\frac{1}{2}\int_{\mathbb{R}^3} |\nabla u_n|^2+\frac{|u|^2}{|x'|^2} dx
 -\frac{t^{6-2s_1}}{6-2s_1}\int_{\mathbb{R}^3}\frac{|u|^{6-2s_1}}{|x|^{s_1}}dx.
 \end{eqnarray*}
 By \eqref{4.5}, we have that as $n\in \mathbb{N}$ large,
 \begin{eqnarray*}
\bar{c}_S+1+\|u_n\|\geqslant\mathcal{J}_S(u_n)-\frac{1}{6-2s_1}\langle\mathcal{J}_S'(u_n),u_n\rangle
=\frac{2-s_1}{6-2s_1}\|u_n\|^2,
\end{eqnarray*}
which is impossible. We have
\begin{eqnarray}\label{4.7}
\begin{array}{ll}
& \displaystyle \frac{2-s_1}{6-2s_1}\|u_n\|^2+\frac{1}{6-2s_1}\langle\mathcal{J}_S'(u_n),u_n\rangle \\
=& \displaystyle \mathcal{J}_S(u_n)+\frac{\lambda(s_2-s_1)}{2(3-s_2)(3-s_1)}\Psi_S(u_n) \int_{\mathbb{R}^3}\frac{|u_n|^{6-2s_2}}{|x|^{s_2}}dx\\[1em] &\displaystyle -\frac{\lambda\|u_n\|^2 }{2(3-s_1)(3-s_2)S^2}\psi'\left(\frac{\|u_n\|^2}{S^2}\right)\int_{\mathbb{R}^3} \frac{|u_n|^{6-2s_2}}{|x|^{s_2}}dx \\
\leqslant &  \displaystyle \mathcal{J}_S(u_n)+\frac{\lambda(s_2-s_1)}{2(3-s_2)(3-s_1)}\Psi_S(u_n) \int_{\mathbb{R}^3}\frac{|u_n|^{6-2s_2}}{|x|^{s_2}}dx\\[1em]
&\displaystyle -\frac{\lambda\|u_n\|^2 }{2(3-s_1)(3-s_2)S^2}\left|\psi'\left(\frac{\|u_n\|^2}{S^2}\right)\right|\int_{\mathbb{R}^3} \frac{|u_n|^{6-2s_2}}{|x|^{s_2}}dx.
\end{array}
\end{eqnarray}
 Suppose, up to a subsequence, that
\begin{eqnarray}
\lim_{n\rightarrow\infty}\|u_n\|\geqslant \frac{S}{2}.\label{4.8}
 \end{eqnarray}
Since $\{u_n\}$ is bounded, it follows from \eqref{4.5} and  \eqref{4.8} that for $n$ sufficiently large,
\begin{eqnarray}\label{4.9}
\frac{2-s_1}{6-2s_1}\|u_n\|^2+\frac{1}{6-2s_1}\langle\mathcal{J}_S'(u_n),u_n\rangle\geqslant CS^2-\frac{\|\mathcal{J}_S'(u_n)\|\|u_n\|}{6-2s_1}\geqslant CS^2-S.
\end{eqnarray}
Notice that if $\|u_n\|>S$ then $\Psi_S(u_n)=0$ and $\psi'\left(\frac{\|u_n\|^2}{S^2}\right)=0.$ So we can obtain that as $n$ large enough, using \eqref{4.7},
\begin{eqnarray*}
CS^2-S\leq \bar{c}_S+1,
\end{eqnarray*}
which is impossible as $S$ large enough. Now we consider the case $\|u_n\|\leq S$. Using $\psi\in C_0^{\infty}(\mathbb{R},[0,1])$, we deduce that
\begin{eqnarray}\label{4.10}
\Psi_S(u_n)\int_{\mathbb{R}^3}\frac{|u_n|^{6-2s_2}}{|x|^{s_2}} dx\leqslant C\|u_n\|^{6-2s_2}\Psi_S(u_n)\leqslant C S^{6-2s_2},
\end{eqnarray}
\begin{eqnarray}\label{4.11}
\frac{\|u_n\|^2}{S^2}\left|\psi'\left(\frac{\|u_n\|^2}{S^2}\right)\right|\int_{\mathbb{R}^3}\frac{|u_n|^{6-2s_2}}{|x|^{s_2}}dx\leqslant C\frac{\|u_n\|^{8-2s_2}}{S^2}\left|\psi'\left(\frac{\|u_n\|^2}{S^2}\right)\right|\leqslant C S^{6-2s_2}.
\end{eqnarray}
By the definitions  of $\bar{c}_S$ and $\tilde{v}$, we infer that
\begin{eqnarray}\label{4.12}
\left.
  \begin{array}{ll}
\bar{c}_S&\displaystyle \leqslant\max_{t\in[0,1]}\mathcal{J}_S(t\tilde{v})\\
&\displaystyle \leqslant\max_{t\in[0,1]}\left\{\frac{t^2}{2}\|\tilde{v}\|^2-\frac{t^{6-2s_1}}{6-2s_1}\int_{\mathbb{R}^3}\frac{|\tilde{v}|^{6-2s_1}}{|x|^{s_1}}dx\right\}\\
&\displaystyle \ \ \ \ +\max_{t\in[0,1]}\left\{\frac{-\lambda t^{6-2s_2}}{6-2s_2}\Psi_S(t\tilde{v}) \int_{\mathbb{R}^3}\frac{|\tilde{v}|^{6-2s_2}}{|x|^{s_2}} dx\right\}.
  \end{array}
\right.
\end{eqnarray}
As in \eqref{4.10}, we derive that
\begin{eqnarray*}
\max_{t\in[0,1]}\left\{\frac{-\lambda t^{6-2s_2}}{6-2s_2}\Psi_S(t\tilde{v})\int_{\mathbb{R}^N}\frac{|\tilde{v}|^{6-2s_2}}{|x|^{s_2}} dx\right\}\leqslant\max_{t\in[0,1]}\left\{\frac{-\lambda C t^{6-2s_2}}{6-2s_2}\|\tilde{v}\|^{6-2s_2}\Psi_S(t\tilde{v})\right\}\leqslant -\lambda C  S^{6-2s_2}.
\end{eqnarray*}
Combining this with \eqref{4.12}, we deduce that $$\bar{c}_S\leqslant C-C\lambda S^{6-2s_2}.$$
This together with $\mathcal{J}_S(u_n)\rightarrow\bar{c}_S$ as $n\rightarrow\infty$ imply that
for $n$ large enough,  \begin{eqnarray}\label{4.13}
 \mathcal{J}_S(u_n)\leqslant C-C\lambda S^{6-2s_2}.
  \end{eqnarray}
Substituting  \eqref{4.10}--\eqref{4.13} in \eqref{4.7}, we have that for $n$ sufficiently large,
 \begin{eqnarray}\label{4.14}
\frac{2-s_2}{6-2s_2}\|u_n\|^2+\frac{1}{6-2s_2}\mathcal{J}_S'(u_n)[u_n]\leqslant C-C\lambda S^{6-2s_2}.
\end{eqnarray}
From \eqref{4.9} and \eqref{4.14}, we obtain that
\begin{eqnarray}\label{4.15}
C-C\lambda S^{6-2s_2}\geqslant CS^2-S ,
\end{eqnarray}
where $C>0$ is independent of  $S$ and $\lambda$.  The inequality \eqref{4.15} would not hold for $S>0$ sufficiently large and $0>\lambda>-S^{2s_2-6}$. The proof is complete. $\hfill\Box$

\begin{lem}\label{lem4.3}
There exists $\lambda^{**}<0$ such that as $\lambda\in(\lambda^{**},0)$, the functional $\I_\lambda$ has a bounded Palais-Smale sequence $\{u_n\}$ at the level $\bar{c}_S$.
\end{lem}
\pf  By Lemma \ref{lem4.1}, the cut-off functional $\mathcal{J}_S$ has a mountain pass level $\bar{c}_S>0$ given by \eqref{4.4} and
 there exists $\{u_n\}\subset \X$ satisfying \eqref{4.5} for each fixed $S>0$ and $\lambda<0$.
 According to Lemma \ref{lem4.2}, we choose $S>0$ large enough and $\lambda^{**}=\lambda^{**}(S)<0$ large such that
 for any $\lambda<\lambda^{**}$,
 $$\limsup_{n\rightarrow\infty}\|u_n\|<\frac{S}{2}.$$
 Combining this with \eqref{4.2} and the definition of $\Psi_{S}$ given in \eqref{4.3}, we derive that for $n$ large enough, $ \mathcal{J}_{S}(u_n)=\mathcal{I}(u_n)$ and $\mathcal{J}'_{S}(u_n)=\mathcal{I}'(u_n).$ Therefore we have
 $ \mathcal{I} (u_n) \rightarrow \bar{c}_S>0$ and $\mathcal{I} '(u_n)\rightarrow 0$ as $ n\rightarrow\infty.$   The proof is complete.
     $\hfill\Box$

By the same proofs of Lemmas \ref{lem3.8}, \ref{lem3.9}, \ref{lem3.10}. We can obtain the following lemmas, in which we investigate the properties of the $\mathrm{(PS)}_{\bar{c}_S}$ sequence $\{u_n\}$ of $\I_\lambda$ found in Lemma \ref{lem4.3}.

 \begin{lem}\label{lem4443}
 If $u_n\rightharpoonup 0$ in $\X$, then for any  domain $B_{a,b}(0)$, up to a subsequence and still denoted by $\{u_n\}$ such that
\begin{eqnarray}
\left.
  \begin{array}{ll} \displaystyle
\int_{B_{a,b}(0)}|\nabla u_n|^2+\frac{|u_n|^2}{|x'|^2}dx\to 0, \ \ \int_{B_{a,b}(0)}\frac{|u_n|^{6-2s_i}}{|x|^{s_i}}dx\to 0, \ \ \ i=1,2.
  \end{array}
\right.\label{4.16}
\end{eqnarray}
\end{lem}
\pf For any $R>r>0$, the compactness of the embedding
\begin{eqnarray}
\X\hookrightarrow\hookrightarrow L^{6-2s_1}(B_{r,R}(0);|x|^{-s_1})\label{4.17}
\end{eqnarray}
implies that
\begin{eqnarray*}
\int_{B_{a,b}(0)}\frac{|u_n|^{6-2s_1}}{|x|^{s_1}}dx \to 0  \ \ \ {\rm as}\ n\to\infty.
\end{eqnarray*}
Applying the function $\eta$ defined in Lemma \ref{lem3.4}.  Since $\eta^2u_n\in \X$
for all $n\in \mathbb{N}$, combining with \eqref{4.17} and \eqref{3.10}, we get that
\begin{eqnarray}
\left.
  \begin{array}{ll}
o(1)&=\langle\I_\lambda'(u_n),\eta^2u_n\rangle\\
&\displaystyle =\int_{\mathbb{R}^3}\nabla u_n\cdot \nabla(\eta^2u_n)+\frac{|\eta u_n|^2}{|x'|^2}dx
-\int_{\mathbb{R}^3}\frac{|u_n|^{6-2s_1}\eta^2}{|x|^{s_1}}dx
-\lambda\int_{\mathbb{R}^3}\frac{|u_n|^{6-2s_2}\eta^2}{|x|^{s_2}}dx\\
&= \displaystyle \int_{\mathbb{R}^3}\eta^2|\nabla u_n|^2+\frac{|\eta u_n|^2}{|x'|^2}dx-\lambda\int_{\mathbb{R}^3}\frac{|u_n|^{6-2s_2}\eta^2}{|x|^{s_2}}dx.
  \end{array}
\right.\label{4.18}
\end{eqnarray}
Since $\lambda<0$, we get, as $n\to\infty$,
\begin{eqnarray*}
\int_{\mathbb{R}^3}|\eta\nabla u_n|^2+\frac{|\eta u_n|^2}{|x'|^2}dx\to0.
\end{eqnarray*}
The  definition $\eta|_{B_{a,b}(0)}\equiv1$ implies
\begin{eqnarray*}
\int_{B_{a,b}(0)}|\nabla u_n|^2+\frac{|u_n|^2}{|x'|^2}dx\to0\ {\rm as}\ n\to\infty.
\end{eqnarray*}
We complete the proof. \hbx

\begin{lem}\label{lem4.5}
If $u_n\rightharpoonup 0$ in $\X$, then there exists a subsequence{\rm(}still denoted by $\{u_n\}${\rm)} such that
\begin{eqnarray*}
{\rm either} \ \ \lim_{n\to\infty}\int_{B_\delta(0')\times B_\delta(0_3)}\frac{|u_n|^{6-2s_1}}{|x|^{s_1}}dx=0 \ \ \ {\rm or} \ \
\lim_{n\to\infty}\int_{B_\delta(0')\times B_\delta(0_3)}\frac{|u_n|^{6-2s_1}}{|x|^{s_1}}dx\geqslant S_{s_1}^{\frac{3-s_1}{2-s_1}}
\end{eqnarray*}
for all $\delta>0$.
\end{lem}

Now, by the transform \eqref{3.19}, then we find that  $\{\tilde{u}_n\}\subset\X$ is also a ${\rm(PS)_{\bar{c}_S}}$ sequence of $\I_\lambda$.

\begin{lem}\label{lem4.6}  There exists $\epsilon_1\in (0,S_{s_1}^{\frac{3-s_1}{2-s_1}}]$ such that for all $\epsilon\in(0, \frac{\epsilon_1}{2})$, there exists  a sequence $\{r_n>0\}$
 such that  $\{\tilde{u}_n\}$  verifies
\begin{eqnarray}
\int_{B_1(0')\times B_1(0_3)}\frac{|\tilde{u}_n|^{6-2s_1}}{|x|^{s_1}}dx=\epsilon. \label{4.20}
\end{eqnarray}
\end{lem}
\pf Set
\begin{eqnarray*}
\mathcal{A}(u_n) \to \mathcal{A}_\infty, \ \  \mathcal{B}(u_n)\to \mathcal{B}_\infty, \ \ \mathcal{C}(u_n)\to \mathcal{C}_\infty.
\end{eqnarray*}
We claim that  $B_\infty>0$, otherwise, since $\{u_n\}$ is a ${\rm (PS)_{\bar{c}_S}}$ sequence, we obtain
 \begin{eqnarray}
\mathcal{A}_\infty-\mathcal{B}_\infty-\lambda \mathcal{C}_\infty=0.\label{4.21}
\end{eqnarray}
Furthermore,
 \begin{eqnarray*}
\mathcal{A}_\infty=\mathcal{B}_\infty+\lambda \mathcal{C}_\infty\leq \mathcal{B}_\infty.
\end{eqnarray*}
Assume $B_\infty=0$, then $A_\infty=0$ and follows from \eqref{4.21}, there exists  a contradiction with the fact $\bar{c}_S>0$. The remaining proof is similar with the proof of Lemma \ref{lem3.10}.\hbx

Now we are ready to complete the proof of Theorem \ref{thm1.4}.

\noindent{\bf Proof of Theorem \ref{thm1.4}.} \ By Lemma \ref{lem4.3},  $\I_\lambda$ has a  ${\rm(PS)_{\bar{c}_S}}$ sequence $\{u_n\}\subset \X$. By Lemma \ref{lem4.6},  the sequence $\{\tilde{u}_n\}$ define by \eqref{3.19}  is also  a bounded ${\rm(PS)_{\bar{c}_S}}$ sequence of $\I_\lambda$ in $\X$.
 Thus there is $\tilde{u}_0\in \X$ such that
\begin{eqnarray*}
  \left\{\begin{array}{ll}
\tilde{u}_n\rightharpoonup \tilde{u}_0 \ \ &\mbox{in} \ \X;\\
\tilde{u}_n\rightharpoonup \tilde{u}_0 \ \ &\mbox{in} \ L^{6-2s_i}(\mathbb{R}^3;|x|^{-s_i}), \ i=1, 2; \\
\tilde{u}_n(x)\rightarrow \tilde{u}_0(x) \ \ & \mbox{a.e. \ on} \  \mathbb{R}^3.
  \end{array}
\right.
\end{eqnarray*}
It follows that $\tilde{u}_0$ is a solution of \eqref{1.17}.  Let $v_n:=\tilde{u}_n-\tilde{u}_0$.
Then $\{v_n\}$ is bounded in $\X$. Define
\begin{eqnarray*}
\mathcal{A}(v_n) \to \mathcal{A}_\infty, \ \  \mathcal{B}(v_n)\to \mathcal{B}_\infty, \ \ \mathcal{C}(v_n)\to \mathcal{C}_\infty.
\end{eqnarray*}
Then by Brezis-Lieb Lemma \cite{1983B-Lieb}, we get
\begin{eqnarray}
 && \I_\lambda(v_n)\to \frac{1}{2}\mathcal{A}_\infty-\frac{1}{6-2s_1}\mathcal{B}_\infty
-\frac{\lambda}{6-2s_2}\mathcal{C}_\infty=\bar{c}_S-\I_\lambda(\tilde{u}_0),\label{4.22} \\
&&\langle\I_\lambda'(v_n),v_n\rangle\to \mathcal{A}_\infty-\mathcal{B}_\infty-\lambda\mathcal{C}_\infty=0.\label{4.23}
\end{eqnarray}
 If $\mathcal{A}_\infty=0$, then $\I_\lambda(\tilde{u}_0)=\hat{c}_S$ and $\tilde{u}_0$ is a nontrivial solution of \eqref{1.17}. Assume $\mathcal{A}_\infty>0$ and
 $\tilde{u}_0=0$. Then  Lemma \ref{lem4.5} implies that
\begin{eqnarray*}
  {\rm  either}\  \displaystyle\lim_{n\to\infty}\int_{B_1(0)}\frac{|\tilde{u}_n|^{6-2s_1}}{|x|^{s_1}}dx=0\ {\rm or}\  \displaystyle\lim_{n\to\infty}\int_{B_1(0)}\frac{|\tilde{u}_n|^{6-2s_1}}{|x|^{s_1}}dx \geqslant S_{s_1}^{\frac{3-s_1}{2-s_1}}.
\end{eqnarray*}
By Lemma \ref{lem4.6}, this is a contradiction to \eqref{4.20} as $0<\epsilon<\frac{1}{2}S_{s_1}^{\frac{3-s_1}{2-s_1}}$. It must be $\tilde{u}_0\not=0$ and $\tilde{u}_0$ is a nontrivial solution of \eqref{1.17}. Combining with Remark \ref{rem}, the equation \eqref{E} has a nontrivial solution. \hfill $\Box$

\section{Asymptotic behavior}
In this section we focus on the proofs of Theorem \ref{thm1.6} in subsection \ref{sec5.1}, Theorem \ref{thm1.7} in subsection \ref{sec5.2}, Theorem \ref{thm1.8} in subsection \ref{sec5.3}, Theorem \ref{thm1.9} in subsection \ref{sec5.4}.

According to Theorems \ref{thm1.1}, \ref{thm1.2}, \ref{thm1.3}, we know that
\begin{eqnarray*}
&&m_{\lambda}=\inf_{u\in \mathcal{S}_{\lambda}} \I_{\lambda}
\end{eqnarray*}
with
\begin{eqnarray*}
&&\mathcal{S}_{\lambda}:=\left\{u\in \X| \I_{\lambda}'(u)=0\right\}
\end{eqnarray*}
is well defined under the conditions of Theorems \ref{thm1.1}, \ref{thm1.2}, \ref{thm1.3}.

\subsection{The case of $0<s_1<2, s_2=2$}\label{sec5.1}

We now prove  Theorem \ref{thm1.6}. We first give a key lemma about the asymptotic behavior of energy level.
\begin{lem}\label{lem5.1}
Assume the assumptions of Theorem \ref{thm1.6} hold. Then there exists a subsequence $\{\lambda_n\}$ satisfying $\lim_{n\to\infty}\lambda_n=0$, such that
\begin{eqnarray*}
\lim_{n\to\infty}m_{\lambda_n}=m_0:=\frac{2-s_1}{6-2s_1} S_{s_1}^{\frac{3-s_1}{2-s_1}}.
\end{eqnarray*}
\end{lem}
\pf We know that there exists a $u_0\in \X$ satisfies $m_0=\I_0(u_0)$ and $u_0$ belongs to mountain pass type solution of $\I_0$. And $m_0=\bar{c}_0=\hat{c}_0$, where $\bar{c}_0=\hat{c}_0$ is defined in \eqref{3.1} and \eqref{3.2}.

Now we prove this lemma in two different situations: $\lambda_n\to 0^+, \lambda_n\to 0^-$.

We first consider the case $\lambda_n>0$ and $\lambda_n\to 0^+$ as $n\to\infty$. Combining with \eqref{3.4} and \eqref{3.5}, we get
\begin{eqnarray}
m_0>m_{\lambda_n}.\label{5.1}
\end{eqnarray}
It follows from Theorem \ref{thm1.1} that, for any $0<\lambda_n<\bar{\lambda}$,  there exists a solution $u_n$ satisfies $m_{\lambda_n}=\I_{\lambda_n}(u_n)$ and $u_n$  is mountain pass type, by \eqref{5.1}, It is not difficult to verify that $u_n$ is bounded in $\X$. Now we prove the fact that $\lim_{n\to\infty}m_{\lambda_n}\geq m_0$. Since, as $n$ large,
\begin{eqnarray*}
0&\geq&\langle\I_{\lambda_n}'(u_n),u_n\rangle\\
&=&\intR|\nabla u_n|^2+\frac{|u_n|^2}{|x'|^2}dx-\intR\frac{|u_n|^{6-2s_1}}{|x|^{s_1}}dx-\lambda_n\intR\frac{|u_n|^{2}}{|x|^{2}}dx\\
&\geq& \|u_n\|^2-S_{s_1}^{s_1-3}\|u_n\|^{6-2s_1}-\frac{\lambda_n}{\bar{\lambda}}\|u_n\|^2\\
&=&\left(1-\frac{\lambda_n}{\bar{\lambda}}\right) \|u_n\|^2-S_{s_1}^{s_1-3}\|u_n\|^{6-2s_1}.
\end{eqnarray*}
It follows from $\lambda_n<\bar{\lambda}$ that there exists $M:=M(s_1)$ such that $\|u_n\|\geq M$. Furthermore,
\begin{eqnarray}
\left.
  \begin{array}{ll}
o(1)+\I_{\lambda_n}(u_n)&=\I_{\lambda_n}(u_n)-\frac{1}{6-2s_1}\langle\I_{\lambda_n}'(u_n),u_n\rangle\\
&=\frac{2-s_1}{6-2s_1}\left(\intR|\nabla u_n|^2+\frac{|u_n|^2}{|x'|^2}dx-\lambda_n\intR\frac{|u_n|^2}{|x|^2}dx\right)\\
&\geq \frac{2-s_1}{12-4s_1}M^2
  \end{array}
\right.\label{5.2}
\end{eqnarray}
as $n\to\infty$. The facts \eqref{5.1} and \eqref{5.2} imply that $\I_{\lambda_n}(u_n)$ is bounded. Thus there exists a subsequence(still denoted by origin mark) such that
\begin{eqnarray*}
m_{\lambda_n}=\I_{\lambda_n}(u_n)\to c>0, \I_{\lambda_n}'(u_n)=0.
\end{eqnarray*}
If $c\geq m_0$, then the proof is complete. Otherwise, $c<m_0$, we will construct a contradiction. The boundedness of sequence $\{u_n\}$ implies that there hold,
\begin{eqnarray*}
\left\{
  \begin{array}{ll}
u_n\rightharpoonup u_0\ {\rm in}\ \X,\\
u_n\rightharpoonup u_0\ {\rm in}\ L^{6-2s_1}(\RT;|x|^{-s_1}),\\
u_n\rightharpoonup u_0\ {\rm in}\ L^{2}(\RT;|x|^{-2}),\\
u_n(x)\to u_0(x)\ {\rm a.e.\ on}\ \RT.
  \end{array}
\right.
\end{eqnarray*}
If follows that $u_0$ is a critical point of $\I_0$ and $\I_0(u_0)\geq0$. Let $v_n=u_n-u_0$, applying the Brezis-Lieb lemma, we can get that
\begin{eqnarray*}
&&\I_{\lambda_n}(v_n)\to c-\I_0(u_0),\\
&&\|v_n\|^2-\intR \frac{|v_n|^{6-2s_1}}{|x|^{s_1}}dx-\lambda_n \intR \frac{|v_n|^2}{|x|^2}dx\to0,\\
&&\lambda_n \intR \frac{|v_n|^2}{|x|^2}dx\to0.
\end{eqnarray*}
We may therefore assume that
\begin{eqnarray*}
\|v_n\|^2\to b,\ \intR \frac{|v_n|^{6-2s_1}}{|x|^{s_1}}dx\to b.
\end{eqnarray*}
The inequality \eqref{Inequality} implies that $b\geq S_{s_1} b^{\frac{1}{3-s_1}}$, which leads to that either $b=0$ or $b\geq S_{s_1}^{\frac{3-s_1}{2-s_1}}$. The case $b=0$ implies that $u_0$ is a nontrivial solution and $\I_0(u_0)=c\geq m_0$, which is a contradiction. However, if $b\geq S_{s_1}^{\frac{3-s_1}{2-s_1}}$, we get that
\begin{eqnarray*}
\displaystyle c\geq \lim_{n\to\infty}\I_{\lambda_n}(v_n)\geq\frac{2-s_1}{6-2s_1}S_{s_1}^{\frac{3-s_1}{2-s_1}}=m_0,
\end{eqnarray*}
which is a contradiction with $c<m_0$. Thus $c\geq m_0$ and the proof of the case $\lambda_n\to 0^+$ is over.

Now we consider the case $\lambda_n<0$ and $\lambda_n\to0^-$ as $n\to\infty$. For any $\lambda<0$, it follows from Theorem \ref{thm1.1} that there exists a   ground state solution $u_\lambda$ of problem \eqref{1.17}, then we have that
\begin{eqnarray}
m_\lambda=\I_\lambda(u_\lambda)\geq \I_\lambda(t_0 u_\lambda)\geq \I_0(t_0 u_\lambda)=\frac{2-s_1}{6-2s_1}\frac{\left(\intR|\nabla u_\lambda|^2+\frac{|u_\lambda|^2}{|x'|^2}dx\right)^{\frac{3-s_1}{2-s_1}}}{\left(\intR\frac{|u_\lambda|^{6-2s_1}}{|x|^{s_1}}dx
\right)^{\frac{1}{2-s_1}}}\geq m_0,\label{5.3}
\end{eqnarray}
where
\begin{eqnarray*}
t_0=\left(\frac{\intR|\nabla u_\lambda|^2+\frac{|u_\lambda|^2}{|x'|^2}dx}{\intR\frac{|u_\lambda|^{6-2s_1}}{|x|^{s_1}}dx}\right)^{\frac{1}{2-s_1}}.
\end{eqnarray*}
For $\lambda_1<\lambda_2\leq 0$, let $u_{\lambda_1}$ be a ground state solution, then
\begin{eqnarray}
m_{\lambda_1}=\I_{\lambda_1}(u_{\lambda_1})\geq \I_{\lambda_1}(t_{\lambda_2}u_{\lambda_1})\geq \I_{\lambda_2}(t_{\lambda_2}u_{\lambda_1})\geq m_{\lambda_2},\label{5.4}
\end{eqnarray}
where $t_{\lambda_2}$ satisfies that $t_{\lambda_2}u_{\lambda_1}\in \mathcal{N} $. Thus for $-1<\lambda<0$, we get that $m_{\lambda}<m_{-1}$. Furthermore, combining with \eqref{5.3} and \eqref{5.4}, one has $\lim_{\lambda\to0^-}m_\lambda=m_0$.

Therefore, we get that there exists a subsequence $\{\lambda_n\}$ satisfying the conclusion of lemma. \hbx

{\bf Proof of Theorem \ref{thm1.6}.} Let $\{u_n\}$ be a ground state solution of problem \eqref{1.17} with $\lambda=\lambda_n$. Then
\begin{eqnarray}
\langle \I_{\lambda_n}'(u_n),u_n\rangle =0\ \ {\rm and} \ \ \I_{\lambda_n}(u_n)\to m_0,\label{5.5}
\end{eqnarray}
where using  Lemma \ref{lem5.1}. Since, for $n$ large,
\begin{eqnarray}
\left.
  \begin{array}{ll}
m_0+1&\geq \I_{\lambda_n}(u_n)-\frac{1}{6-2s_1}\langle \I_{\lambda_n}'(u_n),u_n\rangle\\
&=\frac{2-s_1}{6-2s_1}\left(\|u_n\|^2-\lambda_n \intR\frac{|u_n|^2}{|x|^2}dx\right)\\
&\geq \frac{2-s_1}{6-2s_1}\left(1-\frac{\lambda_n}{\bar{\lambda}}\right)\|u_n\|^2.
  \end{array}
\right.\label{5.6}
\end{eqnarray}
It follows from $\lambda_n<\bar{\lambda}$ that $\{u_n\}$ is bounded in $\X$. Thus there exists a $u_0\in\X$ such that
\begin{eqnarray*}
\left\{
  \begin{array}{ll}
u_n\rightharpoonup u_0\ {\rm in}\ \X,\\
u_n\rightharpoonup u_0\ {\rm in}\ L^{6-2s_1}(\RT;|x|^{-s_1}),\\
u_n\rightharpoonup u_0\ {\rm in}\ L^2(\RT;|x|^{-2}),\\
u_n(x)\to u_0(x)\ {\rm a.e.\ on}\ \RT.
  \end{array}
\right.
\end{eqnarray*}
It follows that $u_0$ is a critical point of $\I_0$ and $\I_0(u_0)\geq0$.

If $u_0\not=0$, we have $\I_0(u_0)\geq m_0$, set $v_n=u_n-u_0$, applying the Brezis-Lieb lemma, we can get that
\begin{eqnarray}
&&\I_{\lambda_n}(v_n)\to m_0-\I_0(u_0),\label{5.7}\\
&&\|v_n\|^2-\intR \frac{|v_n|^{6-2s_1}}{|x|^{s_1}}dx-\lambda_n \intR \frac{|v_n|^2}{|x|^2}dx\to0. \label{5.8}
\end{eqnarray}
Similar with \eqref{5.6}, combining with \eqref{5.7} and \eqref{5.8}, we get that $\I_0(u_0)\leq m_0$. Thus $\I_0(u_0)=m_0$. We may therefore assume that
\begin{eqnarray*}
\|v_n\|^2\to b,\ \intR \frac{|v_n|^{6-2s_1}}{|x|^{s_1}}dx\to b.
\end{eqnarray*}
Therefore,
\begin{eqnarray*}
o(1)&=&\I_{\lambda_n}(v_n)\\
&=&\I_{\lambda_n}(v_n)-\frac{1}{6-2s_1}\langle \I_{\lambda_n}'(u_n),u_n\rangle\\
&=&\frac{2-s_1}{6-2s_1}\left(\|v_n\|^2-\lambda_n\intR \frac{|v_n|^2}{|x|^2}dx\right)\\
&\to& \frac{2-s_1}{6-2s_1}b\ \ \ \  {\rm as}\ n\to\infty.
\end{eqnarray*}
Thus $b=0$, and $u_n\to u_0$, \eqref{5.5} implies that $\I_0(u_0)=m_0$, that is $u_0$ is a ground state solution of equation \eqref{1.17} with $\lambda=0$.

If $u_0=0$, since $m_0>0$, there exists $\epsilon_1\in (0,\frac{1}{2}S_{s_1}^{\frac{2-s_1}{3-s_1}}]$ such that for all $\epsilon\in(0, \epsilon_1)$, there exists  a sequence $\{r_n>0\}$
 such that  $\{\tilde{u}_n:=r_n^{\frac{1}{2}}u_n(r_nx)\}$  verifies
\begin{eqnarray}
\int_{B_1(0')\times B_1(\tilde{0})}\frac{|\tilde{u}_n|^{6-2s_1}}{|x|^{s_1}}dx=\epsilon \label{5.9}
\end{eqnarray}
and \begin{eqnarray*}
\I_{\lambda_n}(\tilde{u}_n)=\I_{\lambda_n}(u_n)\to m_0,\ \  \I_{\lambda_n}'(\tilde{u}_n)=0\ {\rm in}\ (\X)^*.
\end{eqnarray*}
Moreover, there exists $\tilde{u}_0\in\X$ such that
\begin{eqnarray*}
\left\{
  \begin{array}{ll}
\tilde{u}_n\rightharpoonup \tilde{u}_0\ {\rm in}\ \X,\\
\tilde{u}_n\rightharpoonup \tilde{u}_0\ {\rm in}\ L^{6-2s_1}(\RT;|x|^{-s_1}),\\
\tilde{u}_n\rightharpoonup \tilde{u}_0\ {\rm in}\ L^2(\RT;|x|^{-2}),\\
\tilde{u}_n(x)\to \tilde{u}_0(x)\ {\rm a.e.\ on}\ \RT.
  \end{array}
\right.
\end{eqnarray*}
It follows that $\tilde{u}_0$ is a critical point of $\I_0$ and $\I_0(\tilde{u}_0)\geq0$. Set $\tilde{v}_n=\tilde{u}_n-\tilde{u}_0$, applying the Brezis-Lieb lemma, we can get that
\begin{eqnarray}
&&\I_{\lambda_n}(\tilde{v}_n)\to m_0-\I_0(\tilde{u}_0),\nonumber\\
&&\|\tilde{v}_n\|^2-\intR \frac{|\tilde{v}_n|^{6-2s_1}}{|x|^{s_1}}dx-\lambda_n \intR \frac{|\tilde{v}_n|^2}{|x|^2}dx\to0.\label{5.10}
\end{eqnarray}
We may therefore assume that
\begin{eqnarray*}
\|\tilde{v}_n\|^2\to b,\ \intR \frac{|\tilde{v}_n|^{6-2s_1}}{|x|^{s_1}}dx\to b.
\end{eqnarray*}
The inequality \eqref{Inequality} implies that $b\geq S_{s_1} b^{\frac{1}{3-s_1}}$, which leads to that either $b=0$ or $b\geq S_{s_1}^{\frac{3-s_1}{2-s_1}}$. The case $b=0$ implies that $\tilde{u}_0$ is a nontrivial solution, which is desired. If $b\geq S_{s_1}^{\frac{3-s_1}{2-s_1}}$ and $\tilde{u}_0=0$. Coping the proof of Lemma \ref{lem3.4} with $s_1>0$, we can get that for any domain $B_{a,b}(0)$ and any $b>a>0$, up to a subsequence and still denoted by $\{\tilde{u}_n\}$ such that
\begin{eqnarray}
\left.
  \begin{array}{ll} \displaystyle
\int_{B_{a,b}(0)}|\nabla \tilde{u}_n|^2+\frac{|u_n|^2}{|x'|^{2}}dx\to 0, \ \ \int_{B_{a,b}(0)}\frac{|\tilde{u}_n|^{6-2s_i}}{|x|^{s_i}}dx\to 0, \ \ \ i=1,2.
  \end{array}
\right.\label{5.11}
\end{eqnarray}
Set
\begin{eqnarray}
\left.
  \begin{array}{ll}
 &\displaystyle \tilde{\kappa}_1:=\limsup_{n\to\infty}\int_{B_\delta(0')\times B_\delta(0_3)}\frac{|\tilde{u}_n|^{6-2s_1}}{|x|^{s_1}}dx,\\
&\displaystyle \tilde{\kappa} :=\limsup_{n\to\infty}\int_{B_\delta(0')\times B_\delta(0_3)}|\nabla \tilde{u}_n|^2+\frac{|u_n|^2}{|x'|^{2}}dx.
  \end{array}
\right.\label{5.12}
\end{eqnarray}
Based on \eqref{5.11}, similar with   Lemma \ref{lem3.5}, for any $\delta>0$, we get that
\begin{eqnarray}
\tilde{\kappa}=\tilde{\kappa}_1 \label{5.13}
\end{eqnarray}
and
\begin{eqnarray}
\tilde{\kappa}_1^{\frac{1}{3-s_1}}\leq S_{s_1}^{-1}\tilde{\kappa}.\label{5.14}
\end{eqnarray}
Combining with \eqref{5.13} and \eqref{5.14}, we have
\begin{eqnarray}
\tilde{\kappa}_1^{\frac{1}{3-s_1}}\leq S_{s_1}^{-1}\tilde{\kappa}_1.\label{5.15}
\end{eqnarray}
Furthermore, we can obtain that
\begin{eqnarray}
{\rm either}\ \displaystyle\lim_{n\to\infty}\int_{B_1(0')\times B_1(\tilde{0})}\frac{|\tilde{u}_n|^{6-2s_1}}{|x|^{s_1}}dx=0\ {\rm or}\ \displaystyle\lim_{n\to\infty}\int_{B_1(0')\times B_1(\tilde{0})}\frac{|\tilde{u}_n|^{6-2s_1}}{|x|^{s_1}}dx \geqslant S_{s_1}^{\frac{2-s_1}{3-s_1}}.\label{5.16}
\end{eqnarray}
This is a contradiction with \eqref{5.9} as $0<\epsilon<\frac{1}{2}S_{s_1}^{\frac{2-s_1}{3-s_1}}$. It must be $\tilde{u}_0\not=0$. Thus $\I_0(\tilde{u}_0)\geq m_0$. And $\lim_{n\to\infty}\I_{\lambda_n}(\tilde{v}_n)\leq 0$, combining with \eqref{5.10},
\begin{eqnarray*}
o(1)&\geq &\I_{\lambda_n}(\tilde{v}_n)\\
&=&\I_{\lambda_n}(\tilde{v}_n)-\frac{1}{6-2s_1}\langle \I_{\lambda_n}'(\tilde{v}_n),\tilde{v}_n\rangle\\
&=&\frac{2-s_1}{6-2s_1}\|\tilde{v}_n\|^2-\frac{\lambda_n(2-s_1)}{6-2s_1}\intR \frac{|\tilde{v}_n|^2}{|x|^2}dx\\
&\to& \frac{2-s_1}{6-2s_1}b,
\end{eqnarray*}
which is a contradiction with $b\geq S_{s_1}^{\frac{3-s_1}{2-s_1}}$. Thus $\tilde{u}_n\to\tilde{u}_0$ in $\X$.  The proof of Theorem \ref{thm1.6} is over combining with Remark \ref{rem}.  \hbx

\subsection{The case of $\lambda>0, 0<s_1<s_2<2$}\label{sec5.2}

We now prove  Theorem \ref{thm1.7}. We  consider the asymptotic behavior of energy level as follows.
\begin{lem}\label{lem5.2}
Assume $\lambda>0$ holds. Then there exists a subsequence $\{\lambda_n\}$ satisfying $\lim_{n\to\infty}\lambda_n=0$, such that
\begin{eqnarray*}
\lim_{n\to\infty}m_{\lambda_n}=m_0:=\frac{2-s_1}{6-2s_1} S_{s_1}^{\frac{3-s_1}{2-s_1}}.
\end{eqnarray*}
\end{lem}
\pf We know that there exists a $u_0\in \X$ satisfies $m_0=\I_0(u_0)$ and $u_0$ belongs to mountain pass type solution of $\I_0$. And $m_0=\bar{c}_0=\hat{c}_0$, where $\bar{c}_0, \hat{c}_0$ are defined in \eqref{3.1} and \eqref{3.2}.
Combining with \eqref{3.4} and \eqref{3.5}, we get
\begin{eqnarray}
m_0>m_{\lambda_n}.\label{5.17}
\end{eqnarray}
It follows from Theorem \ref{thm1.2} that, for any $\lambda_n>0$,  there exists a solution $u_n$ satisfies $m_{\lambda_n}=\I_{\lambda_n}(u_n)$ and $u_n$  is mountain pass type, by \eqref{5.17}, It is not difficult to verify that $u_n$ is bounded in $\X$. Now we prove the fact that $\lim_{n\to\infty}m_{\lambda_n}\geq m_0$. Since, as $n$ large,
\begin{eqnarray*}
1&\geq&\langle\I_{\lambda_n}'(u_n),u_n\rangle\\
&=&\intR|\nabla u_n|^2+\frac{|u_n|^2}{|x'|^2}dx-\intR\frac{|u_n|^{6-2s_1}}{|x|^{s_1}}dx-\lambda_n\intR\frac{|u_n|^{6-2s_2}}{|x|^{s_2}}dx\\
&\geq& \|u_n\|^2-S_{s_1}^{s_1-3}\|u_n\|^{6-2s_1}-S_{s_2}^{s_2-3}\|u_n\|^{6-2s_2}.
\end{eqnarray*}
It follows from that there exists $M:=M(s_1,s_2)$ such that $\|u_n\|\geq M$. Furthermore,
\begin{eqnarray}
\left.
  \begin{array}{ll}
o(1)+\I_{\lambda_n}(u_n)&=\I_{\lambda_n}(u_n)-\frac{1}{6-2s_2}\langle\I_{\lambda_n}'(u_n),u_n\rangle\\
&=\frac{2-s_2}{6-2s_2}\|u_n\|^2+\frac{s_2-s_1}{2(3-s_1)(3-s_2)}\intR\frac{|u_n|^{2^*(s_1)}}{|x|^{s_1}}dx\\
&\geq \frac{2-s_2}{12-4s_2}M^2
  \end{array}
\right.\label{5.18}
\end{eqnarray}
as $n\to\infty$. The facts \eqref{5.17} and \eqref{5.18} imply that $\I_{\lambda_n}(u_n)$ is bounded. Thus there exists a subsequence(still denoted by origin mark) such that
\begin{eqnarray*}
m_{\lambda_n}=\I_{\lambda_n}(u_n)\to c>0, \I_{\lambda_n}'(u_n)=0.
\end{eqnarray*}
If $c\geq m_0$, then the proof is complete. Otherwise, $c<m_0$, we will construct a contradiction. The boundedness of sequence $\{u_n\}$ implies that there hold,
\begin{eqnarray*}
\left\{
  \begin{array}{ll}
u_n\rightharpoonup u_0\ {\rm in}\ \X,\\
u_n\rightharpoonup u_0\ {\rm in}\ L^{6-2s_i}(\RT;|x|^{-s_i}),\ i=1,2,\\
u_n(x)\to u_0(x)\ {\rm a.e.\ on}\ \RT.
  \end{array}
\right.
\end{eqnarray*}
If follows that $u_0$ is a critical point of $\I_0$ and $\I_0(u_0)\geq0$. Let $v_n=u_n-u_0$, applying the Brezis-Lieb lemma, we can get that
\begin{eqnarray*}
&&\I_{\lambda_n}(v_n)\to c-\I_0(u_0),\\
&&\|v_n\|^2-\intR \frac{|v_n|^{6-2s_1}}{|x|^{s_1}}dx-\lambda_n \intR \frac{|v_n|^{6-2s_2}}{|x|^{s_2}}dx\to0,\\
&&\lambda_n \intR \frac{|v_n|^{6-2s_2}}{|x|^{s_2}}dx\to0
\end{eqnarray*}
We may therefore assume that
\begin{eqnarray*}
\|v_n\|^2\to b,\ \intR \frac{|v_n|^{6-2s_1}}{|x|^{s_1}}dx\to b.
\end{eqnarray*}
The inequality \eqref{Inequality} implies that $b\geq S_{s_1} b^{\frac{1}{3-s_1}}$, which leads to that either $b=0$ or $b\geq S_{s_1}^{\frac{3-s_1}{2-s_1}}$. The case $b=0$ implies that $u_0$ is a nontrivial solution and $\I_0(u_0)=c\geq m_0$, which is a contradiction. However, if $b\geq S_{s_1}^{\frac{3-s_1}{2-s_1}}$, we get that
\begin{eqnarray*}
\displaystyle c\geq \lim_{n\to\infty}\I_{\lambda_n}(v_n)\geq\frac{2-s_1}{6-2s_1}S_{s_1}^{\frac{3-s_1}{2-s_1}}=m_0,
\end{eqnarray*}
which is a contradiction with $c<m_0$. Thus $c\geq m_0$ and the proof is over.\hbx

{\bf Proof of Theorem \ref{thm1.7}.} Let $\{u_n\}$ be a ground state solution of problem \eqref{1.17} with $\lambda=\lambda_n>0$. Then
\begin{eqnarray}
\langle \I_{\lambda_n}'(u_n),u_n\rangle =0\ {\rm and}\ \I_{\lambda_n}(u_n)\to m_0,\label{5.19}
\end{eqnarray}
where using  Lemma \ref{lem5.2}. Since, for $n$ large,
\begin{eqnarray}
\left.
  \begin{array}{ll}
m_0+1&\geq \I_{\lambda_n}(u_n)-\frac{1}{6-2s_2}\langle \I_{\lambda_n}'(u_n),u_n\rangle\\
&=\frac{2-s_2}{6-2s_2}\|u_n\|^2+\frac{s_2-s_1}{2(3-s_1)(3-s_2)}\intR \frac{|u_n|^{2^*(s_1)}}{|x|^{s_1}}dx.
  \end{array}
\right.\label{5.20}
\end{eqnarray}
It follows that $\{u_n\}$ is bounded in $\X$. Thus there exists a $u_0\in\X$ such that
\begin{eqnarray*}
\left\{
  \begin{array}{ll}
u_n\rightharpoonup u_0\ {\rm in}\ \X,\\
u_n\rightharpoonup u_0\ {\rm in}\ L^{6-2s_i}(\RT;|x|^{-s_i}),\ i=1,2,\\
u_n(x)\to u_0(x)\ {\rm a.e.\ on}\ \RT.
  \end{array}
\right.
\end{eqnarray*}
It follows that $u_0$ is a critical point of $\I_0$ and $\I_0(u_0)\geq0$.

If $u_0\not=0$, we have $\I_0(u_0)\geq m_0$, set $v_n=u_n-u_0$, applying the Brezis-Lieb lemma, we can get that
\begin{eqnarray}
&&\I_{\lambda_n}(v_n)\to m_0-\I_0(u_0),\label{5.21}\\
&&\|v_n\|^2-\intR \frac{|v_n|^{6-2s_1}}{|x|^{s_1}}dx-\lambda_n \intR \frac{|v_n|^{6-2s_2}}{|x|^{s_2}}dx\to0. \label{5.22}
\end{eqnarray}
Similar \eqref{5.20}, combining with \eqref{5.21} and \eqref{5.22}, we get that $\I_0(u_0)\leq m_0$. Thus $\I_0(u_0)=m_0$. We may therefore assume that
\begin{eqnarray*}
\|v_n\|^2\to b,\ \intR \frac{|v_n|^{6-2s_1}}{|x|^{s_1}}dx\to b.
\end{eqnarray*}
Therefore,
\begin{eqnarray*}
o(1)&=&\I_{\lambda_n}(v_n)\\
&=&\I_{\lambda_n}(v_n)-\frac{1}{6-2s_1}\langle \I_{\lambda_n}'(u_n),u_n\rangle\\
&=&\frac{2-s_1}{6-2s_1}\|v_n\|^2+\frac{\lambda_n(s_1-s_2)}{2(3-s_1)(3-s_2)}\intR \frac{|v_n|^{2^*(s_2)}}{|x|^{s_2}}dx\\
&\to& \frac{2-s_1}{6-2s_1}b\ \ \ \  {\rm as}\ n\to\infty.
\end{eqnarray*}
Thus $b=0$ and $u_n\to u_0$, \eqref{5.19} implies that $\I(u_0)=m_0$, that is $u_0$ is a ground state solution of equation \eqref{1.17} with $\lambda=0$.

If $u_0=0$, since $m_0>0$, there exists $\epsilon_1\in (0,\frac{1}{2}S_{s_1}^{\frac{3-s_1}{2-s_1}}]$ such that for all $\epsilon\in(0, \epsilon_1)$, there exists  a sequence $\{r_n>0\}$
 such that  $\{\tilde{u}_n:=r_n^{\frac{1}{2}}u_n(r_nx)\}$  verifies
\begin{eqnarray}
\int_{B_1(0')\times B_1(\tilde{0})}\frac{|\tilde{u}_n|^{6-2s_1}}{|x|^{s_1}}dx=\epsilon \label{5.23}
\end{eqnarray}
and \begin{eqnarray*}
\I_{\lambda_n}(\tilde{u}_n)=\I_{\lambda_n}(u_n)\to m_0, \I_{\lambda_n}'(\tilde{u}_n)=0\ {\rm in}\ (\X)^*.
\end{eqnarray*}
Moreover, there exists $\tilde{u}_0\in\X$ such that
\begin{eqnarray*}
\left\{
  \begin{array}{ll}
\tilde{u}_n\rightharpoonup \tilde{u}_0\ {\rm in}\ \X,\\
\tilde{u}_n\rightharpoonup \tilde{u}_0\ {\rm in}\ L^{6-2s_1}(\RT;|x|^{-s_1}),\\
\tilde{u}_n(x)\to \tilde{u}_0(x)\ {\rm a.e.\ on}\ \RT.
  \end{array}
\right.
\end{eqnarray*}
It follows that $\tilde{u}_0$ is a critical point of $\I_0$ and $\I_0(\tilde{u}_0)\geq0$. Set $\tilde{v}_n=\tilde{u}_n-\tilde{u}_0$, applying the Brezis-Lieb lemma, we can get that
\begin{eqnarray}
&&\I_{\lambda_n}(\tilde{v}_n)\to m_0-\I_0(\tilde{u}_0),\nonumber\\
&&\|\tilde{v}_n\|^2-\intR \frac{|\tilde{v}_n|^{6-2s_1}}{|x|^{s_1}}dx-\lambda_n \intR \frac{|\tilde{v}_n|^{6-2s_2}}{|x|^{s_2}}dx\to0.\label{5.24}
\end{eqnarray}
We may therefore assume that
\begin{eqnarray*}
\|\tilde{v}_n\|^2\to b,\ \intR \frac{|\tilde{v}_n|^{6-2s_1}}{|x|^{s_1}}dx\to b.
\end{eqnarray*}
The inequality \eqref{Inequality} implies that $b\geq S_{s_1} b^{\frac{1}{3-s_1}}$, which leads to that either $b=0$ or $b\geq S_{s_1}^{\frac{3-s_1}{2-s_1}}$. The case $b=0$ implies that $\tilde{u}_0$ is a nontrivial solution, which is desired. If $b\geq S_{s_1}^{\frac{3-s_1}{2-s_1}}$ and $\tilde{u}_0=0$. Similar with \eqref{5.12}-\eqref{5.16} in the proof of Theorem \ref{thm1.6}, we can obtain that
\begin{eqnarray*}
{\rm either}\ \displaystyle\lim_{n\to\infty}\int_{B_1(0')\times B_1(\tilde{0})}\frac{|\tilde{u}_n|^{6-2s_1}}{|x|^{s_1}}dx=0\ {\rm or}\ \displaystyle\lim_{n\to\infty}\int_{B_1(0')\times B_1(\tilde{0})}\frac{|\tilde{u}_n|^{6-2s_1}}{|x|^{s_1}}dx \geqslant S_{s_1}^{\frac{3-s_1}{2-s_1}}.
\end{eqnarray*}
This contradicts to \eqref{5.23} as $0<\epsilon<\frac{1}{2}S_{s_1}^{\frac{3-s_1}{2-s_1}}$. It must be $\tilde{u}_0\not=0$. Thus $\I_0(\tilde{u}_0)\geq m_0$. And $\lim_{n\to\infty}\I_{\lambda_n}(\tilde{v}_n)\leq 0$, combining with \eqref{5.24},
\begin{eqnarray*}
o(1)&\geq &\I_{\lambda_n}(\tilde{v}_n)\\
&=&\I_{\lambda_n}(\tilde{v}_n)-\frac{1}{6-2s_1}\langle \I_{\lambda_n}'(\tilde{v}_n),\tilde{v}_n\rangle\\
&=&\frac{2-s_1}{6-2s_1}\|\tilde{v}_n\|^2+\frac{\lambda_n(s_1-s_2)}{2(3-s_1)(3-s_2)}\intR \frac{|\tilde{v}_n|^{6-2s_2}}{|x|^{s_2}}dx\\
&\to& \frac{2-s_1}{6-2s_1}b,
\end{eqnarray*}
which is a contradiction with $b\geq S_{s_1}^{\frac{3-s_1}{2-s_1}}$. Thus $\tilde{u}_n\to\tilde{u}_0$ in $\X$. Finally, according to Remark \ref{rem} we can prove the Theorem \ref{thm1.7}. \hbx

\subsection{The case of $\lambda<0, 0<s_1<s_2<2$} \label{sec5.3}
In this section we focus on the proof of Theorem \ref{thm1.8}, that is  the case $\lambda<0, s_1<s_2$, we first give a asymptotic behavior of energy level.
\begin{lem}\label{lem5.3}
Assume $ \lambda<0, 0<s_1<s_2<2$. Then there exists sequence $\{\lambda_n\}$ satisfying $\lim_{n\to\infty}\lambda_n=0$ such that
\begin{eqnarray*}
\lim_{n\to\infty}m_{\lambda_n}=m_0:=\frac{2-s_1}{6-2s_1} S_{s_1}^{\frac{3-s_1}{2-s_1}}.
\end{eqnarray*}
\end{lem}
\pf We know that there exists a $u_0\in\X$ satisfies $m_0=\I_0(u_0)$ and $u_0$ belongs to mountain pass type solution of $\I_0$. And $m_0=\bar{c}_0=\hat{c}_0$, where $\bar{c}_0, \hat{c}_0$ are defined in \eqref{3.1} and \eqref{3.2}.
For any $\lambda<0$, it follows from Theorem \ref{thm1.3} that there exists a mountain pass type ground state solution $u_\lambda$ of problem \eqref{1.17}, then we have that
\begin{eqnarray}
m_\lambda=\I_\lambda(u_\lambda)\geq \I_\lambda(t_0 u_\lambda)\geq \I_0(t_0 u_\lambda)=\frac{2-s_1}{6-2s_1}\frac{\left(\intR|\nabla u_\lambda|^2+\frac{|u_\lambda|^2}{|x'|^2}dx\right)^{\frac{3-s_1}{2-s_1}}}{\left(\intR\frac{|u_\lambda|^{6-2s_1}}{|x|^{s_1}}dx
\right)^{\frac{1}{2-s_1}}}\geq m_0,\label{5.25}
\end{eqnarray}
where
\begin{eqnarray*}
t_0=\left(\frac{\intR|\nabla u_\lambda|^2+\frac{|u_\lambda|^2}{|x'|^2}dx}{\intR\frac{|u_\lambda|^{6-2s_1}}{|x|^{s_1}}dx}\right)^{\frac{1}{2-s_1}}.
\end{eqnarray*}
For $\lambda_1<\lambda_2\leq 0$, let $u_{\lambda_1}$ be a ground state solution, then
\begin{eqnarray}
m_{\lambda_1}=\I_{\lambda_1}(u_{\lambda_1})\geq \I_{\lambda_1}(t_{\lambda_2}u_{\lambda_1})\geq \I_{\lambda_2}(t_{\lambda_2}u_{\lambda_1})\geq m_{\lambda_2},\label{5.26}
\end{eqnarray}
where $t_{\lambda_2}$ satisfies that $t_{\lambda_2}u_{\lambda_1}\in \mathcal{N} $. Thus for $-1<\lambda<0$, we get that $m_{\lambda}<m_{-1}$. Thus, combining with \eqref{5.25} and \eqref{5.26}, one has $\lim_{\lambda\to0^-}m_\lambda=m_0$.\hbx

{\bf Proofs of Theorem \ref{thm1.8}.} The proof is similar with Theorem \ref{thm1.7}. \hbx

\subsection{Proof of Theorem \ref{thm1.9}}\label{sec5.4}
In the present section, we only prove Theorem \ref{thm1.9}. According to  Theorem \ref{thm1.4}, we get that for large $S>0$, there exists a $\lambda^*<0$ such that as $\lambda^*<\lambda<0$, the equation \eqref{1.17} has  a nontrivial solution  $\tilde{u}_\lambda$ satisfying
\begin{eqnarray*}
\J_S(\tilde{u}_\lambda)=\I_\lambda(\tilde{u}_\lambda)=\hat{c}_\lambda,
\end{eqnarray*}
where
\begin{eqnarray*}
  \hat{c}_\lambda=\inf_{\gamma\in \Gamma_S^\lambda} \max_{t\in [0, 1]} \J_S(\gamma(t))>0
  \end{eqnarray*}
and
  $$\Gamma_S^\lambda =\left\{\gamma\in C([0, 1], \X): \ \gamma(0)=0, \ \J_S(\gamma(1))<0\right\}.$$
Fixed $S$ large, if $\lambda^1<\lambda^2$ and $\gamma\in \Gamma_S^{\lambda^1}$, it follows from
\begin{eqnarray*}
\J_{S}(\gamma(1))=\frac{1}{2}\mathcal{A}(\gamma(1))-\frac{1}{6-2s_1}\mathcal{B}(\gamma(1))-\frac{\lambda^1}{6-2s_2}
\Psi_{S}(\gamma(1))\mathcal{C}(\gamma(1))<0
\end{eqnarray*}
that
\begin{eqnarray*}
\J_{S}(\gamma(1))=\frac{1}{2}\mathcal{A}(\gamma(1))-\frac{1}{6-2s_1}\mathcal{B}(\gamma(1))-\frac{\lambda^2}{6-2s_2}
\Psi_{S}(\gamma(1))\mathcal{C}(\gamma(1))<0.
\end{eqnarray*}
Thus  $\Gamma_S^{\lambda^1}\subset \Gamma_S^{\lambda^2}$ and so $\hat{c}_{\lambda^2}\leq \hat{c}_{\lambda^1}$. Thus there exists a sequence $\{\lambda_n>0\}$ satisfying $\lim_{n\to+\infty}\lambda_n=0$  such that $\lim_{n\to+\infty}\hat{c}_{\lambda_n}=\hat{c}_0>0$.
We use  $\tilde{u}_{\lambda_n}$ to denote the solutions of \eqref{1.17} corresponding to the energy $\hat c_{\lambda_n}$, that is
\begin{eqnarray*}
\I_{\lambda_n}(\tilde{u}_{\lambda_n})=\hat{c}_{\lambda_n}, \ \  \I_{\lambda_n}'(\tilde{u}_{\lambda_n})=0.
\end{eqnarray*}
Coping the proof of Theorem \ref{thm1.1}, we see that $\{u_n=r_n^{\frac{1}{2}}\tilde{u}_{\lambda_n}(r_nx)\}$ also a solution of  \eqref{1.17} and there exists a $u\in \X\setminus\{0\}$ such that $u_n\to u$ in $\X$. As a consequence, $u$ satisfies   \eqref{1.17} with $\lambda=0$ and the energy $\hat{c}_0$. Finally, according to Remark \ref{rem},  the Theorem \ref{thm1.9} is proved. \hbx

\section*{Acknowledgment}
This work is supported by NSFC(12301258,12271373,12171326) and Sichuan Science and Technology Program(2023NSFSC1298).

\section*{Data availibility}
Data sharing not applicable to this article as no datasets were generated or analysed during the current study.

\section*{Declarations}
\subsection*{Conflict of interest} The authors have no Conflict of interest to declare that are relevant to the content of this
article.

\end{document}